\documentclass[a4paper,12pt]{article}
\usepackage{amsmath,amsthm,amssymb,fancybox,longtable,graphics,amscd,multirow}
\usepackage[dvipdfm]{graphicx}
\usepackage{caption,setspace,lscape,hyperref,color,array}
\usepackage[all]{xy}

\numberwithin{equation}{section}

\newtheorem{thm}{\textsc{Theorem}}[section]
\newtheorem{prop}{\textsc{Proposition}}[section]
\newtheorem{lem}{\textsc{Lemma}}[section]
\newtheorem{cor}{\textsc{Corollary}}[section]
\newtheorem{PropDefn}{\textsc{Proposition-Definition}}[section]

\newtheorem{Introthm}{\textsc{Theorem}}

\theoremstyle{definition}
\newtheorem{eg}{Example \vspace{2mm}}[section]
\newtheorem{remark}{\textsc{Remark}}
\newtheorem{defn}{\textsc{Definition}}[section]

\def\proof{\textsc{Proof. }}
\def\ackn{{\bf Acknowledgements.\,}}
\def\QED{$\Box$}
\def\END{$\blacksquare$}
\def\mbi#1{\boldsymbol{#1}} 

\def\Pic{\mathop{\mathrm{Pic}}\nolimits}

\def\sign{\sigma}

\def\Div{\mathop{\mathrm{div}}\nolimits}

\def\f1#1{\frac{1}{#1}}

\begin{document}
\setlength{\baselineskip}{16pt}

\title{A note on simple $K3$ singularities and families of weighted $K3$ surfaces}
\author{Makiko Mase \\ \small{mtmase@arion.ocn.ne.jp} }
\date{\empty}
\maketitle

\noindent
\begin{small}
https://orcid.org/0000-0002-8927-0714 
\end{small}

\begin{abstract}
We discuss properties of the Seifert form for simple $K3$ singularities, and of the Picard lattices of families of weighted $K3$ surfaces. 
We study a collection $\mathcal{M}_{(\rho,\,\delta)}$ of $K3$ surfaces polarized by their Picard lattices that are in the set $\mathcal{L}_{(\rho,\,\delta)}$ of certain lattices. 
We also report a numerical formula that relates the Seifert form for the singularities and the Picard lattice of the family. 
\end{abstract}


\noindent
\ackn
The author appreciates to whom it may concern to my assistantship at the University of Mannheim for their hospitality and for giving her an opportunity of precious experiences while studying the subjects. 

\section{Introduction}
Being classified by Yonemura~\cite{Yo90}, certain simple $K3$ singularities are isolated quasi-homogeneous hypersurface singularities $(f=0)$ defined by quasi-homogeneous polynomials in four variables.  
The weighted projective spaces with weights that are given in the classification have anticanonical divisors that produce families of {\it weighted} $K3$ surfaces, in the sense that the minimal models of general sections are $K3$ surfaces that are polarized by the Picard lattice of generic sections of the family. 
In~\cite[Corollary 4.4, Remark 4.5]{Yo90}, it is remarked that if one takes a defining polynomial of a simple $K3$ singularity to be general, there is a relation between the rank of the resolution graph in normal $K3$ surface and the number of parameters of the moduli of $K3$ surfaces. 
It is well-known by Torelli-type thoerem that the Picard lattice, with the intersection form on the singular second cohomology group, plays an important role in the surjectivity of the period map and to describe the moduli space of lattice-polarized $K3$ surfaces. 

On the other hand, the Seifert form is defined on the Milnor lattice, the singular second homology group of $f^{-1}(1)$, of the singularity. 
Steenbrink~\cite{S77-2} studies many important topological invariants concerning isolated hypersurface singularities such as the Jacobi algebra and the Hodge numbers of the singularity, and in particular, gives an explicit formula of the Poincar\'e series of a quasi-homogeneous IHS. 
Saeki~\cite{Sae00} studies the $\mathbb{R}$- and $\mathbb{C}$-extensions of Seifert form of isolated quasi-homogeneous singularities to conclude those extensions determine a topology of them. 

The aim of the article is to study the collection $\mathcal{M}_{\rho,\,\delta}$ of $K3$ surfaces $S=S_L$ polarized by a lattice $L\in\mathcal{L}_{(\rho,\,\delta)}$ such that $\Pic{(S)} = L$, where, $\mathcal{L}_{(\rho,\,\delta)}$ is the set of non-degenerate even hyperbolic lattices that are primitive sublattices of the $K3$ lattice $\Lambda_{K3}$, of signature $(1,\,\rho-1)$, and discriminant number $\delta$ for all pairs $(\rho,\,\delta)$ that appear as the rank and discriminant number of each Picard lattice of the family of weighted $K3$ surfaces associated to simple $K3$ singularities, given by Belcastro~\cite{Be02}. 
Our first main theorem is the following: 
\begin{Introthm}[Theorems \ref{LatticeSimple} and \ref{MainThmLattice}]
\begin{itemize}
\item[$(a)$] There exists a unique lattice with signature $(1,\,\rho-1)$ and of discriminant number $\delta$ up to isomorphism, if a pair $(\rho,\,\delta)$ is given with 
   $(1)$ $\delta$ is a prime number, as in Table~\ref{Prime}; 
   $(2)$ $\rho=1$, as in Table~\ref{rho1}; and 
   $(3)$ $\delta=1$, as in Table~\ref{delta1}. 
\item[$(b)$] For a given pair $(\rho,\,\delta)$ as in Table~\ref{Reclassify}, one can give all quadratic forms of signature $2-\rho$ and defined on the discriminant group of order $\delta$ up to isomorphism. 
Moreover, for each triple $(\rho,\,\delta,\, q)$, there exists a unique lattice of signature $(1,\,\rho-1)$ that admits $q$ as its discriminant quadratic form. 
\end{itemize}
\end{Introthm}
\noindent
And we conclude: 

\noindent
{\bf Corollary~\ref{AllElliptic}}
{\it Any $K3$ surfaces in $\mathcal{M}_{\rho,\,\delta}$ for all $(\rho,\,\delta)$ with $\rho\not=1,2$ have a structure of elliptic fibration. }\\

We also consider a relation between Seifert form and Picard lattice associated to simple $K3$ singularities as ``lattice'' structures rather than the famous Poincar\'e duality. 
More precisely, we will report a numerical formula which is stated as follows: 
Let $c_2$ be the dimension of the eigenspace of eigenvalue $1$ of the real Seifert form of a simple $K3$ singularity, and $\rho$ the Picard number of the family of $K3$ surfaces corresponding to the singularity. 
Denote by $l(\Delta^{[1]})$ the total number of lattice points on all edges of the Newton polytope of the defining polynomial of the general section in the family. \\

\noindent
{\bf Proposition~\ref{MainThmSeifert}}
{\it The following formula holds for all simple $K3$ hypersurface singularities. }
\[
c_2 = l(\Delta^{[1]})-3 = 20-\rho.  
\]

We explain the structure of the article. 
In Section~\ref{Setup}, we recall the concepts of the Picard lattice as well as weighted $K3$ surfaces together with an introduction to some basic definitions and facts. 

In Section~\ref{StudyOfLattices}, we first classify lattices in $\mathcal{L}_{(\rho,\,\delta)}$ by a standard lattice theory recalled in the first part, then we give our conclusion concerning members in $\mathcal{M}_{(\rho,\,\delta)}$ in Corollary~\ref{AllElliptic}. 

In Section~\ref{LatticesSingularities}, after recalling the formulas for quasi-homogeneous singularities of the Poincar\'e series due to Steenbrink, and of the real Seifert form by Saeki, following an example, we prove our result Proposition~\ref{MainThmSeifert}. 
\\

We close the introduction by emphasizing the following remark: 
\begin{remark}[On a comparison of Belcastro's and our results]
Among Table~\ref{Reclassify}, there are both the same and different entries of~\cite{Be02}. 

\cite{Be02} starts from the weights $(a_0,\,a_1,\,a_2,\,a_3)$ classified by~\cite{Yo90} and then, computed the lattices. 
Thus, the elliptic fibrations in~\cite{Be02}'s list are those that are determined by the structure of the weighted projective space with weights $(a_0,\,a_1,\,a_2,\,a_3)$. 

On the other hand, in our article, we start from the pairs $(\rho,\,\delta)$ of the rank $\rho$ and the discriminant number $\delta$ picked up in~\cite{Be02}. 
So, it is natural to obtain the same entries as in~\cite{Be02} in our Table~\ref{Reclassify}. 
However, at the same time, there are some $\mathbb{Z}$-modules with given $(\rho,\,\delta)$ but different quadratic forms (equivalently, different intersection matrix). 
In other words, our elliptic fibres are those that are not necessarily determined by the structure of the weighted projective space with weights $(a_0,\,a_1,\,a_2,\,a_3)$ of~\cite{Yo90} . 
This is the difference from~\cite{Be02}. 
\end{remark}

\section{The Picard lattice of a $K3$ surface and family of weighted $K3$ surfaces}\label{Setup}
\begin{defn}[The anticanonical divisor]
Let $X$ be a normal projective complex algebraic variety of dimension $n$. 
Denote by $X_{\rm reg}$ the nonsingular part of $X$, and $j:X_{\rm reg}\hookrightarrow X$ the natural inclusion map. 
Let $\Omega_{X_{\rm reg}}$ be the sheaf of differentiable $1$-forms on $X_{\rm reg}$. 
\begin{itemize}
\item[(i)] The {\it canonical sheaf} $\omega_X$ of $X$ is defined by 
\[
\omega_X := j_*\left(\displaystyle\bigwedge^n\Omega_{X_{\rm reg}}\right). 
\]
\item[(ii)] The {\it canonical divisor} of $X$ is the divisor $K_X$ on $X$ whose associated sheaf is isomorphic to the canonical sheaf $\omega_X$. 
The canonical divisor is unique up to linear equivalence. 
\item[(iii)] The {\it anticanonical divisor} of $X$ is the divisor ${-}K_X$ that is the $(-1)$ times of the canonical divisor. 
\item[(iv)] We denote by $|{-}K_X|$ the linear system associated to the anticanonical divisor, and call the {\it anticanonical linear system}. 
\END
\end{itemize}
\end{defn}

\begin{defn}[$K3$ surface]
Let $S$ be a complex algebraic variety of dimension $2$. 
The variety $S$ is called a {\it $K3$ surface} if the canonical divisor and irregularity of $S$ are both trivial. 
\END
\end{defn}

\begin{defn}[The weighted projective space]
Let $\mbi{a}=(a_0,\,a_1,\,a_2,\,a_3)$ be a quadruple of positive integers. 
Take polynomials $f$ and $g$ from the polynomial ring $\mathbb{C}[z_0,\,z_1,\,z_2,\,z_3]$. 
\begin{itemize}
\item[(i)] Define a relation $\sim$ on $\mathbb{C}[z_0,\,z_1,\,z_2,\,z_3]$ by 
\[
f(z_0,\,z_1,\,z_2,\,z_3) \sim g(z_0,\,z_1,\,z_2,\,z_3) 
\]
if 
\[
 f(z_0,\,z_1,\,z_2,\,z_3) = g(\zeta^{a_0}z_0,\,\zeta^{a_1}z_1,\,\zeta^{a_2}z_2,\,\zeta^{a_3}z_3)
\]
holds for all non-zero complex number $\zeta$. 
\item[(ii)] The {\it weighted projective space} with {\it weights $a_0$,\,$a_1$,\,$a_2$,\,$a_3$} is the projective variety denoted by 
\[
\mathbb{P}(\mbi{a}) = \mathbb{P}(a_0,\,a_1,\,a_2,\,a_3)
\]
that is the projectivization of the polynomial ring $\mathbb{C}[z_0,\,z_1,\,z_2,\,z_3]$ with the relation $\sim$. 
We say that the {\it weight of variable $z_i$ is $a_i$}\, $(i=0,1,2,3)$. 
\item[(iii)]
A non-degenerate polynomial $f(z_0,\,z_1,\,z_2,\,z_3)$ in $\mathbb{C}[z_0,\,z_1,\,z_2,\,z_3]$ with the relation $\sim$ is a {\it quasi-homogeneous polynomial of degree $d$} if 
\[
f(\zeta^{a_0} z_0,\,\zeta^{a_1} z_1,\,\zeta^{a_2} z_2,\,\zeta^{a_3} z_3) = \zeta^df(z_0,\,z_1,\,z_2,\,z_3)
\]
holds for any non-zero complex number $\zeta$. 
We call a tuple $(a_0,\,a_1,\,a_2,\,a_3; d)$ the {\it weight system}. 
\END
\end{itemize}
\end{defn}

\begin{defn}[Well-formed weight system]
The weights $a_0$,\,$a_1$,\,$a_2$,\,$a_3$ of the weighted projective space $\mathbb{P}(a_0,\,a_1,\,a_2,\,a_3)$ is said {\it well-formed} if the following conditions are satisfied: 
\begin{itemize}
\item[(i)] $\gcd{(a_0,\,a_1,\,a_2,\,a_3)} = 1$, 
\item[(ii)] $\gcd{(a_i,\,a_j,\,a_k)} = 1$ for any distinct $i,j,k\in\{ 0,1,2,3\}$, and 
\item[(iii)] $a_0\leq a_1\leq a_2\leq a_3$. 
\END
\end{itemize}
\end{defn}

\begin{defn}[Isolated hypersurface singularity (abbreviated : IHS)]
Let $f(z_0,\,z_1,\,z_2,\,z_3)$ be a quasi-homogeneous polynomial of degree $d$ in the weighted projective space $\mathbb{P}(a_0,\,a_1,\,a_2,\,a_3)$ with weights $a_0$,\,$a_1$,\,$a_2$,\,$a_3$. 
Let $X$ be the set of all solutions of the equation $f=0$. 
Then, the germ $(X,\,0)$ is a singularity. 
We call the singularity $(X,\,0)$ an {\it isolated hypersurface singularity}. 
\END
\end{defn}

We recall some fundamental facts on the lattices associated to $K3$ surfaces. 
Let $S$ be a $K3$ surface. 

\begin{defn}[Lattices]
\begin{itemize}
\item[(i)] A {\it lattice} is a non-degenerate $\mathbb{Z}$-module with a symmetric bilinear form on it. 
\item[(ii)] The {\it dual lattice of a lattice $L$}, denoted by $L^*$, is the holomorphism group $Hom(L,\,\mathbb{Z})$. 
There is a natural embedding of groups $L^*\hookrightarrow L$. 
\item[(iii)] The quotient group $A_L:=L\slash L^*$ is called the {\it discriminant group}. 
The group $A_L$ admits a quadratic form. 
\item[(iv)] A lattice $L$ is {\it unimodular} if the discriminant group is the trivial group. 
\item[(v)] The {\it signature of a lattice} is a pair $(s_{(+)},\, s_{(-)})$ of the number $s_{(+)}$ of positive eigenvalues, and the number $s_{(-)}$ of negative eigenvalues of the bilinear form. 
\item[(vi)] The {\it rank of a lattice} is the number of a generator of the lattice as a $\mathbb{Z}$-module.  
\END
\end{itemize}
\end{defn}
\begin{defn}[The $K3$ lattice, and Picard and transcendental lattices of a $K3$ surface]
Let $L$ and $L'$ be two lattices. 
We denote by $L\oplus L'$ the direct sum of two lattices. 

We denote by $U$ the even unimodular lattice of signature $(1,1)$, and $E_8$ the even unimodular lattice of signature $(0,8)$. 
\begin{itemize}
\item[(i)] The lattice $\Lambda_{K3}:=U^{\oplus 3}\oplus E_8^{\oplus 2}$ is called the {\it $K3$ lattice}. 
\item[(ii)] The cohomology group $H^1(S,\,\mathcal{O}^*_S)$ is called the {\it Picard lattice of $S$}, and is denoted by $\Pic(S)$. 
The rank of the Picard lattice of $S$ is called the {\it Picard number of $S$}. 
\item[(iii)] The orthogonal complement of the Picard lattice in the $K3$ lattice is called the {\it transcendental lattice of $S$}, and denoted by $T(S)$. (See the discussion below. )
\END
\end{itemize}
\end{defn}

We explain how the Picard lattice $\Pic(S)$ of $S$ admits a structure of a lattice inherited from the $K3$ lattice. 

By a natural short exact sequence 
\[
\xymatrix{
0 \ar[r] & \mathbb{Z} \ar[r] & \mathcal{O}_S \ar[r] & \mathcal{O}_S^* \ar[r] & 1
}
\]
we get a long exact sequence of cohomology groups
\[
\xymatrix{
\cdots \ar[r] & H^1(S,\,\mathcal{O}_S) \ar[r] & H^1(S,\, \mathcal{O}_S^*) \ar[r]^c & H^2(S,\,\mathbb{Z})\ar[r] & \cdots. 
}
\]
Since $H^1(S,\,\mathcal{O}_S)=0$, that is, the irregularity is trivial, the map $c$ is injective. 
It is known that $H^2(S,\,\mathbb{Z})$ is an even unimodular lattice of signature $(3,\,19)$, thus by classification of lattices, it is isomorphic to the $K3$ lattice. 
Moreover the Picard lattice is a primitive sublattice of $H^2(S,\,\mathbb{Z})\simeq\Lambda_{K3}$. 
The signature of $\Pic(S)$ is $(1,\,\rho-1)$, where $\rho$ is the Picard number of $S$, and the signature of the transcendental lattice $T(S)$ of $S$ is $(2,\, 20-\rho)$. 

\bigskip

We recall a family of weighted $K3$ surfaces associated to certain simple $K3$ singularities. 
From now on, we use the following notations: \\
$\mbi{w}=(w_0,\ldots,w_n)$ : (rational) weight system, \\
where $w_i := s_i\slash t_i$ with $\gcd{(s_i,\,t_i)}=1$ for $i=0,\ldots,n$, and $N:=\mathrm{lcm}{(t_0,\ldots,t_n)}$. \\
$Q=(q_0,\ldots,q_n)$ with $q_i:=Nw_i$ : (integral) weight system, and $\vert Q\vert := \sum_{i=0}^n q_i$. \\

We employ as the definition of a simple $K3$ singularity the geometric statement due to~\cite{IW92}. 

\begin{defn}[Simple $K3$ singularity]
A {\it simple $K3$ singularity} is a $3$-dimensional normal isolated singularity $(X,\,x)$ satisfying the following conditions: \\
Let $E$ be the exceptional divisor of a $\mathbb{Q}$-factorial minimal model. 
\begin{itemize}
\item[(i)] The divisor $E$ is irreducible, and 
\item[(ii)] the divisor $E$ is a normal surface that is birationally equivalent to a $K3$ surface. 
\end{itemize}
\END
\end{defn}

We can consider simple $K3$ singularities that are isolated hypersurface singularities. 
According to~\cite{Yo90}, they are defined by a polynomial $f=f(z_1,\,z_2,\,z_3,\,z_4)$ of degree $1$ in the graded polynomial ring $\mathbb{C}[z_1,\,z_2,\,z_3,\,z_4]$ with weight of variable $z_i$ being $w_i$ for $i=1,2,3,4$. 
And these weight systems are classified by~\cite{Yo90}. 

It can be translated with an integral weight system. 
We may assume that a quadruple $(q_1,\, q_2,\,q_3,\,q_4)$ of integers is well-formed. 
Thus we can define the weighted projective space $\mathbb{P}(\mbi{w}):=\mathbb{P}(q_1,\, q_2,\,q_3,\,q_4)$ as the projectivization of the graded ring $\mathbb{C}[z_1,\,z_2,\,z_3,\,z_4]$ with weights $q_1$,\,$q_2$,\,$q_3$,\,$q_4$. 
The anticanonical divisor of $\mathbb{P}(\mbi{w})$ is described by~\cite{DolgachevWP} : 
\[
\mathcal{O}_{\mathbb{P}(\mbi{w})}({-}K_{\mathbb{P}(\mbi{w})}) \simeq \mathcal{O}_{\mathbb{P}(\mbi{w})}(q_1+q_2+q_3+q_4). 
\]
In our case, one can verify that $N:=\mathrm{lcm}{(t_0,\ldots,t_n)}$ coincides with the value $q_1+q_2+q_3+q_4$. 
Thus, the quasi-homogeneous polynomials of degree $N$ is parametrized by the anticanonical linear system $|{-}K_{\mathbb{P}(\mbi{w})}|$. 

\begin{defn}
The {\it family of weighted $K3$ surfaces} is the collection of varieties parametrized by $|{-}K_{\mathbb{P}(\mbi{w})}|$, where $\mbi{w}$ is one of the weights in the list by~\cite{Yo90}. 
\END
\end{defn}

It is known that any general member of the family of weighted $K3$ surfaces is birationally equivalent to a $K3$ surface. 
By~\cite{Bruzzo-Grassi}, the Picard lattices of the minimal model of any generic members in a family of weighted $K3$ surfaces are isomorphic, thus, we define the following: 

\begin{defn}
Let $\mathcal{F}$ be a family of weighted $K3$ surfaces. 
The Picard lattice of the minimal model of any generic members in $\mathcal{F}$ is called {\it the Picard lattice of the family $\mathcal{F}$}. 
\END
\end{defn}

\begin{defn}
Let $X$ be a normal projective complex algebraic variety of dimension $3$. 
The variety $X$ is called a {\it Fano $3$-fold} if the anticanonical divisor ${-}K_X$ is ample. 
\END
\end{defn}

In other words, it is equivalent that the weight system $\mbi{w}$ is associated to a simple $K3$ singularity and that the weighted projective space $\mathbb{P}(\mbi{w})$ is a Fano $3$-fold ({\it c.f.},~\cite{ReiC3f}).

\section{A study of lattice-polarized $K3$ surfaces in $\mathcal{M}_{(\rho,\,\delta)}$}\label{StudyOfLattices}
\subsection{Notations and basic facts}
For a non-degenerate even hyperbolic lattice $S$, denote by $\rho=\rho(S)$ the rank of $S$, by $\delta=\delta(S)$ the discriminant number of $S$, by $G=G_S$ the discriminant group of $S$, by $q=q_S$ the quadratic form on $S$. 
We also let $G_q$ be the finite group on which a quadratic form $q$ is defined, and further, $\mathfrak{l}(A)$ the length of a finitely-generated group $A$, that is, the minimal number of generators of $A$. 
``No.'' in tables are the indices of weight systems in~\cite{Yo90}. 
Denote by $L(k)$ with a natural number $k$ and a lattice $L$ admitting a quadratic form $q$, the lattice with quadratic form $kq$. 
We indicate the intersection matrix $M$ with respect to a basis $\{g_1,\ldots,g_\rho\}$ of the lattice of rank $\rho$ by $S=\mathopen\langle g_1,\ldots,g_\rho\mathclose\rangle_{\mathbb{Z}}=(\mathbb{Z}^\rho,\,M)$. 

\begin{figure}[!htb]
\begin{center}
\includegraphics[width=4cm]{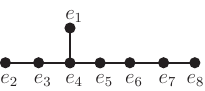}
\begingroup
\captionof{figure}{Lattice of type $E_8$}\label{E8lattice}
\endgroup
\end{center}
\end{figure}

Let $\{ u_1,\,u_2\}$ be a $\mathbb{Z}$-basis of the hyperbolic lattice $U$ of rank $2$ such that $u_i^2=0$ for $i=1,2$ and $u_1.u_2=1$. 
Let $\{ e_1, e_2, e_3, e_4, e_5, e_6, e_7, e_8\}$ be a $\mathbb{Z}$-basis of the negative-definite even unimodular lattice $E_8$ of rank $8$ whose Dynkin diagram is given in Figure~\ref{E8lattice}. 

\subsection{Fundamental quadratic forms}
We refer the reader~\cite[pp.58--59]{Br83}~\cite{Ni80} and~\cite{Be02} for more details. 
We define ``fundamental'' quadratic forms. 

\begin{defn}\label{FundamentalQForms}
The quadratic form $q=w^\varepsilon_{p,k}$ is defined on the group $\mathbb{Z}\slash p^k\mathbb{Z}$ for prime number $p$ and integer $k\geq1$. 

If $p\not=2$, then, $q(1)=a\cdot p^{-k}\mod p\mathbb{Z}$, and $\varepsilon=\left( \frac{a}{p}\right)$ by Jacobi-Legendre symbol. 
That is, $\varepsilon=\pm1$; $\varepsilon=1$ if $a$ is a square, and $\varepsilon=-1$ otherwise.

If $p=2$, then, $q(1)=\varepsilon\cdot2^{-k}\mod2\mathbb{Z}$, and $\varepsilon=\pm1, \pm5$.

The quadratic form $u_k$ (resp. $v_k$) is defined on $(\mathbb{Z}_{2^k})^2$ with intersection matrix
\[
2^{-k}\begin{pmatrix} 0 & 1 \\ 1 & 0 \end{pmatrix}, 
\quad \textrm{resp. }
2^{-k}\begin{pmatrix} 2 & 1 \\ 1 & 2 \end{pmatrix}
\qquad
\textrm{for } k\geq1. \qquad\blacksquare
\]
\end{defn}

By a direct computation, we can show: 

\begin{PropDefn}
The following lattices $S_q$ are respectively representation lattices for the quadratic form $q$. 
\begin{itemize}
\item[$\bullet$] $S_{w^{-5}_{2,3}}  :=  \langle u_1+3u_2+e_6+e_8,\, e_5,\, e_4\rangle_{\mathbb{Z}}, $
\item[$\bullet$] 
$ S_{w^{1}_{2,4}} := \mathopen\Bigg\langle 
\begin{tiny}
\begin{array}{l}
-\sum_{i=4}^7e_i, \, e_7, \, -e_1 -\sum_{i=3}^7e_i, \quad
-\sum_{i\in\{ 1,2,3,6,7,8\}}e_i - 2 \sum_{i=4}^5e_i , \\
-3 \sum_{i=1,6}e_i - 2 \sum_{i=2,7}e_i - 4 \sum_{i=3,5}e_i - 6 e_4  - e_8, \\
\sum_{i=1,8}e_i + 2 \sum_{i=2,6,7}e_i + 3 \sum_{i=3,4,5}e_i, \\
e_8 + 2\sum_{i=2,7}e_i + 4 \sum_{i=1,3,6}e_i  + 6 e_5 + 7 e_4  
\end{array}
\end{tiny}
\mathclose\Bigg\rangle_{\mathbb{Z}}, $
\item[$\bullet$] $S_{w^{5}_{2,4}} := \langle e_7{-}e_8,\, e_6,\, e_5\rangle_{\mathbb{Z}}, $
\item[$\bullet$] $S_{w^{-5}_{2,4}} := \mathopen\Bigg\langle 
\begin{matrix}
e_1, \, e_2, \,  e_3, \, -e_1-e_2-2 e_3-2 e_4-e_5+e_7+e_8, \\ 
-2 e_1-2 e_2-4 e_3-3 e_4-3 e_5-3 e_6-3 e_7-2 e_8
\end{matrix}
\mathclose\Bigg\rangle_{\mathbb{Z}}, $
\item[$\bullet$] $S_{w^{1}_{2,5}} := \mathopen\Bigg\langle 
\begin{array}{l}
\sum_{i=1,2,5,6}e_i+2 \sum_{i=3}^4e_i,\quad \sum_{i=2}^5e_i, \\ 
\sum_{i=1,3,4,5}e_i, \quad -\sum_{i=3,5,6}e_i,\quad \sum_{i=5,6,8}e_i,\\
e_7-e_2, \quad -e_3-e_4-e_8
\end{array}
\mathclose\Bigg\rangle_{\mathbb{Z}}, $
\item[$\bullet$] $S_{w^{5}_{2,5}} :=
\small
\mathopen\Bigg\langle 
\begin{array}{l}
e_8,\quad 2 u_2+ e_1,\quad e_7,\quad e_5,\quad e_6, \\ 
-u_2+\sum_{i=1,2,5} e_i+2 e_3+3 e_4, \\
u_1-4 u_2 -e_3-2 e_4-e_5 \\
\end{array}
\mathclose\Bigg\rangle_{\mathbb{Z}}, $ \quad
\item[$\bullet$] $S_{w^{-5}_{2,5}} := \langle u_1+6u_2+e_1,\, u_2+e_2,\, e_3\rangle_{\mathbb{Z}}, $ 
\item[$\bullet$] $S_{w^{1}_{3,2}} := \langle u_1+u_2,\, 2u_1-u_2 \rangle_{\mathbb{Z}}, $ 
\quad $\bullet$ $S_{w^{-1}_{3,2}} := \langle u_1-u_2,\, u_1+2u_2 \rangle_{\mathbb{Z}}, $ \,
\item[$\bullet$] $S_{w^{1}_{3,3}} := \langle e_2,\, e_1+e_3+e_4-e_5+e_6+e_7-e_8\rangle_{\mathbb{Z}}, $
\item[$\bullet$] $S_{w^{-1}_{3,3}} :=
\small
\mathopen\Bigg\langle 
\begin{tiny}
\begin{array}{l}
-\sum_{i=1,2,5}e_i- 2\sum_{i=3}^4 e_i, \quad
\sum_{i=2,6}e_i + 2\sum_{i=1,3,5} e_i + 3 e_4 , \\
\sum_{i=1,2,7,8}e_i + 2\sum_{i=3,5,6} e_i + 3 e_4 , \quad
-\sum_{i=1,5}^8e_i - 2 e_4, \\
\sum_{i=1,7,8}e_i, \quad
- e_8 -2\sum_{i=1,2,7} e_1 - 3 e_4 - 4\sum_{i=3,5,6} e_i  
\end{array}
\end{tiny}
\mathclose\Bigg\rangle_{\mathbb{Z}}, $
\item[$\bullet$] $S_{w^{-1}_{5,1}} := \langle u_1+u_2,\, u_1+e_4\rangle_{\mathbb{Z}}, $ \quad
$\bullet$ $S_{w^{-1}_{7,1}} := \langle e_1-e_2,\,e_4\rangle_{\mathbb{Z}}, $
\item[$\bullet$] $S_{w^{-1}_{11,1}}  := 
\small
\mathopen\Bigg\langle 
\begin{array}{l}
  e_5 + e_6, \quad  e_7, \quad e_8,\quad e_1 + 2 e_4 + e_5, \\
  e_2 + e_3 - e_6 - e_7 - e_8, \\
  e_1 + e_3 + 2 e_4 + e_5 + 2 e_6 + e_7 
\end{array}
\mathclose\Bigg\rangle_{\mathbb{Z}}, $
\item[$\bullet$] $S_{v_1\perp w^{-1}_{3,1}}  := \mathopen\langle e_1+e_2,\, -e_2-e_5 \mathclose\rangle_{\mathbb{Z}}, $
\item[$\bullet$] $S_{w^{-1}_{19,1}}  := \mathopen\Bigg\langle 
\begin{matrix}
e_5+e_6+e_7, \quad e_8, \quad e_1+e_2+e_3+2 e_4+e_5, \\
e_2+e_3+e_7+e_8, \quad e_3-e_5,e_5-2 e_7-e_8
\end{matrix}
\mathclose\Bigg\rangle_{\mathbb{Z}}. $ \quad \QED
\end{itemize}
\end{PropDefn}

\begin{remark}
Denote by $(t_{(+)},\,t_{(-)})$ the signature of $S_q$. 
\begin{itemize}  
\item[$(1)$] The lattice $S_{w_{7,1}^{-1}}$ is isometric to $(A_6)^\perp_{E_8}$ (see \textnormal{\cite{Nis96}}). 
\item[$(2)$] The lattices $S_{w^{1}_{3,2}}$ and $S_{w^{-1}_{3,2}}$ are sublattices of $U$ of rank 2 with $t_{(+)}=1$. 
\item[$(3)$] The lattices $S_{w^{-5}_{2,3}}$,\, $S_{w^{-5}_{2,5}}$, and $S_{w^{-1}_{5,1}}$ are sublattices of $U\oplus E_8$ with $t_{(+)}=1$. 
\item[$(4)$] The lattices $S_{w^{1}_{2,4}}$,\, $S_{w^{5}_{2,4}}$,\, $S_{w^{-5}_{2,4}}$,\, $S_{w^{1}_{2,5}}$,\, $S_{w^{1}_{3,3}}$,\, $S_{w^{-1}_{3,3}}$,\, $S_{w^{-1}_{7,1}}$,\, $S_{w^{-1}_{11,1}}$,\, $S_{v_1\perp w^{-1}_{3,1}}$, and $S_{w^{-1}_{19,1}}$ are negative-definite sublattices in $E_8$. 
\item[$(5)$] The lattice $S_{w^{5}_{2,5}}$ is a negative-definite sublattice in $U\oplus E_8$. 
\end{itemize}
\end{remark}

For a lattice with quadratic form $q$ of signature $(t_{(+)},\,t_{(-)})$, we can define, in virtue of ~\cite[Theorem 1.3.3]{Ni80}, the {\it signature of the quadratic form $q$} by  
\[
\sign{(q)}:= t_{(+)} - t_{(-)} \mod 8. 
\]
We can define the {\it{orthogonal sum}} of two quadratic forms $q_1,\, q_2$ defined respectively on finite groups $G_1,\, G_2$ by for all $x_1\in G_1,\, x_2\in G_2$, 
\[
(q_1\perp q_2)(x_1,\,x_2) = q_1(x_1) + q_2(x_2). 
\] 
The signature of $q_1\perp q_2$ can be calculated by~\cite[Proposition 1.11.4]{Ni80} as 
\[
\sign{(q_1\perp q_2)} \equiv \sign{(q_1)} + \sign{(q_2)} \mod 8. 
\]

It is fundamental~\cite[$8^\circ$]{Ni80} that any quadratic forms are decomposed into sums of quadratic forms $w^\varepsilon_{p,k},\, u_k$, and $v_k$ and that any lattices are uniquely determined by their rank and quadratic discriminant form. 

\begin{prop}[\cite{Wa64}, {\cite[Proposition 1.11.2]{Ni80}}]\label{SignatureFormula}
The signatures of quadratic forms are given as follows: 
\begin{eqnarray*}
\sign{(w^{\varepsilon}_{p,k})} \equiv k^2(1-p) + 4k\eta\mod 8, & \textnormal{where} & p\not=2, \tiny \left(\frac{\varepsilon}{p}\right)=(-1)^\eta; \\
\sign{(w^{\varepsilon}_{2,k})} \equiv \varepsilon + 4k\omega(\varepsilon)\mod8, &\textnormal{where}& \omega(\varepsilon)\equiv\frac{\varepsilon^2-1}{8}\mod2; \\
\sign{(v_{k})}\equiv4k\mod8; & & \sign{(u_{k})}\equiv0\mod8. \qquad \square
\end{eqnarray*}
\end{prop}

\begin{lem}\label{sgnComp}
\begin{enumerate}
\item[$(1)$] $\sign{(w^{\pm1}_{2,k})}\equiv\pm1\mod8$. \hspace{4mm}$(2)$ $\sign{(w^{\pm5}_{2,2k})}\equiv \pm5\mod8$. 
\item[$(3)$] $\sign{(w^{\pm5}_{2,2k+1})}\equiv\pm1\mod8$. \hspace{5mm}$(4)$ $\sign{(w^{\pm1}_{3,2k})}\equiv 0\mod8$. 
\item[$(5)$] $\sign{(w^{\pm1}_{2,2k+1})}\equiv\mp2\mod8$. 
\end{enumerate}
\end{lem}
\noindent
{\sc Proof. }
We compute the signature with the value $\omega(\varepsilon)$ in Proposition~\ref{SignatureFormula}. 
\begin{enumerate}
\item[(1)] Since $\omega(\pm1)\equiv0\mod2$, we have 
$
\sign{(w^{\pm1}_{2,k})}\equiv\pm1+0\equiv\pm1\mod8. 
$
\item[(2),] (3)\, Since $\omega(\pm5)\equiv1\mod2$, we have \\
$
\sign{(w^{\pm5}_{2,2k})}\equiv\pm5+4\cdot2k\cdot1\equiv\pm5\mod8, 
$\\
$
\sign{(w^{\pm5}_{2,2k+1})}\equiv\pm5+4\cdot(2k+1)\cdot1\equiv\pm5+4=\pm1\mod8. 
$
\item[(4)] We have \quad
$
\sign{(w^{\pm1}_{3,2k})}\equiv(2k)^2(1-3)+4\cdot2k\cdot\eta\mod8 \equiv 0\mod8. 
$
\item[(5)] We have
\begin{eqnarray*}
\sign{(w^{1}_{3,2k+1})} & \equiv & (2k+1)^2(1-3)+4(2k+1)\cdot0\mod8 \\
                           & \equiv & -2\cdot1 \mod8\equiv-2\mod8. \\
\sign{(w^{-1}_{3,2k+1})} & \equiv & (2k+1)^2(1-3)+4(2k+1)\cdot1\mod8 \\
                           & \equiv & -2\cdot1 + 4 \mod8\equiv2\mod8. \qquad \square
\end{eqnarray*}
\end{enumerate}

By Lemma~\ref{sgnComp} and other minor calculations, we can give signatures of quadratic forms in Table~\ref{QuadraticF}, and lattices that may be sublattices of the Picard lattice of a $K3$ surface in Tables \ref{Alattices}, \ref{DElattices}, \ref{hyperbolic} ({\it c.f.}~\cite{Be02}).  
Here, $k$ is a natural number ($\geq1$). 

\renewcommand{\arraystretch}{1.8}

\begin{tiny}
\begin{tabular}{l}
\begin{minipage}{7cm}
\[
\begin{array}{c||p{8mm}|p{8mm}|p{8mm}|p{8mm}|p{8mm}|p{8mm}|p{8mm}|p{8mm}|p{8mm}|p{8mm}|p{8mm}|p{8mm}|p{8mm}|}
q & $w^1_{2,k}$ & $w^{-1}_{2,k}$ & $w^5_{2,2k}$ & $w^{-5}_{2,2k}$ & $w^5_{2,2k+1}$ & $w^{-5}_{2,2k+1}$ & $(w^1_{2,1})^2$ & $(w^{-1}_{2,1})^2$ & $u_{k}$ & $v_{k}$ & $w^1_{3,1}$ & $w^{-1}_{3,1}$  \\
\hline
\sign{(q)} & 1 & 7 & 5 & 3 & 1 & 7 & 2 & 6 & 0 & $4k$ & 6 & 2 
\end{array}
\]
\end{minipage} \\[-10pt]
\begin{minipage}{7cm}
\[
\begin{array}{c||p{8mm}|p{8mm}|p{8mm}|p{8mm}|p{8mm}|p{8mm}|p{8mm}|p{8mm}|p{8mm}|p{8mm}|p{8mm}|p{8mm}|}
q & $w^{\pm1}_{3,2k}$ & $w^1_{3,2k+1}$ & $w^{-1}_{3,2k+1}$ & $w^1_{5,1}$ & $w^{-1}_{5,1}$ & $w^{\pm1}_{5,2}$ & $w^1_{7,1}$ & $w^{-1}_{7,1}$ & $w^1_{11,1}$ & $w^{-1}_{11,1}$ & $w^1_{13,1}$ & $w^{-1}_{13,1}$ \\
\hline
\sign{(q)} & 0 & 6 & 2 & 4 & 0 & 0 & 2 & 6 & 6 & 2 & 4 & 0
\end{array}
\]
\end{minipage} \\[-10pt]
\begin{minipage}{2cm}
\[
\begin{array}{c||p{8mm}|p{8mm}|p{8mm}|p{8mm}|}
q & $w^1_{17,1}$ & $w^{-1}_{17,1}$ & $w^1_{19,1}$ & $w^{-1}_{19,1}$ \\
\hline
\sign{(q)} & 0 & 4 & 6 & 2 
\end{array}
\]
\end{minipage}
\end{tabular}
\begingroup
\captionof{table}{Quadratic forms and signature}\label{QuadraticF}
\endgroup
\vspace{5mm}
\end{tiny} 

\begin{tiny}
\[
\begin{array}{c||c|c|c|c|c|c|c|c|c|c|c|}
S & A_1 &  A_2 & A_3 & A_4 & A_5 & A_6 & A_7 & A_8 & A_9 & A_{10} & A_{11} \\
\hline
\rho(S) & 1 & 2 & 3 & 4 & 5 & 6 & 7 & 8 & 9 & 10 & 11 \\
\hline
\delta(S) & -2 & 3 & -4 & 5 & -6 & 7 & -8 & 9 & -10 & 11 & -12 \\
\hline
G(S) & \mathbb{Z}_2 & \mathbb{Z}_3 & \mathbb{Z}_4 & \mathbb{Z}_5 & \mathbb{Z}_6 & \mathbb{Z}_7 & \mathbb{Z}_8 & \mathbb{Z}_9 & \mathbb{Z}_{10} & \mathbb{Z}_{11} & \mathbb{Z}_{12} \\
\hline
q(S) & w^{-1}_{2,1} & w^1_{3,1} & w^5_{2,2} & w^1_{5,1} & w^1_{2,1}\perp w^{-1}_{3,1} & w^1_{7,1} & w^1_{2,3} & w^1_{3,2} & \mathbb{Z}_{10} & w^1_{11,1} & \mathbb{Z}_{12}
\end{array}
\]
\begingroup
\captionof{table}{Lattices of type $A_n\, (1\leq n\leq 11)$}\label{Alattices}
\endgroup 
\vspace{5mm}
\[
\begin{array}{c||c|c|c|c|c|c|c|c|c|c|c|c|}
S & D_4 & D_5 & D_6 & D_7 & D_8 & D_9 & D_{10} & D_{11} & E_6 & E_7 & E_8 \\
\hline
\rho(S) & 4 & 5 & 6 & 7 & 8 & 9 & 10 & 11 & 6 & 7 & 8 \\
\hline
\delta(S) & 4 & -4 & 4 & -4 & 4 &-4 & 4 & -4 & 3 & -2 & 1 \\
\hline
G(S) & \mathbb{Z}_2^2 & \mathbb{Z}_4 & \mathbb{Z}_2^2 & \mathbb{Z}_4 & \mathbb{Z}_2^2 & \mathbb{Z}_4 & \mathbb{Z}_2^2 & \mathbb{Z}_4 & \mathbb{Z}_3 & \mathbb{Z}_2 & \{ 0\}\\
\hline
q(S) & v_1 & w^{-5}_{2,2} & (w^1_{2,1})^2 & w^1_{2,2} & u_1 & w^{-1}_{2,2} & (w^{-1}_{2,1})^2 & w^5_{2,2} & w^{-1}_{3,1} & w^1_{2,1} & 1
\end{array}
\]
\begingroup
\captionof{table}{Lattices of type $D_m\, (4\leq m\leq 11)$ and type $E_6,\, E_7,\, E_8$}\label{DElattices}
\endgroup 
\vspace{5mm}
\[
\begin{array}{c||c|c|c|c|}
S & U & U(2) & U(3) & U(4) \\
\hline
\delta(S) & -1 & -4 & -9 & -16 \\
\hline
G(S) & \{0\} & \mathbb{Z}_4 & \mathbb{Z}_9 & \mathbb{Z}_{16} \\
\hline
q(S) & 1 & u_1 & w^1_{3,1}\perp w^{-1}_{3,1} & u_2 
\end{array}
\]
\begingroup
\captionof{table}{Hyperbolic lattices of type $U(k),\,1\leq k\leq4$}\label{hyperbolic}
\endgroup 
\end{tiny}

\begin{remark}
Although we fail to find any references, we expect followings by observation. 
\\
${\rm (1)}\quad  q(A_{n-1}) =
\begin{cases}
w^1_{p,1} & \Leftrightarrow p(\not=2) \textnormal{ is a prime number. } \\
w^1_{p,k} & \Leftrightarrow n=p^k, \, p \textnormal{ is a prime number. } \\
\end{cases}
$ \\
\\
{\rm (2)}\, stable equivalences of lattices: $D_m \simeq D_{m+8}$ {\rm($m\geq4$)}, and $A_n\simeq D_{n+8}$ {\rm($n\geq1$)}. 
\end{remark}

For a hyperbolic singularity of type $T_{p,q,r}$, that is, $x^p + y^q + z^r = 0$ with $1/p+1/q+1/r<1$, the lattice $L''_{p,q,r}$, which is of rank $\rho=p+q+r-2$ and of signature $(1,\rho-1,0)$, sits in the isomorphisms 
\[
L''_{p,q,r}\oplus U \simeq L_{p,q,r},\qquad L''_{p,q,r}\oplus K \simeq L'_{p,q,r}, 
\]
where $K$ is a rank-one lattice generated by an isotropy element, and $L'_{p,q,r}$ is the Milnor lattice of the singularity~\cite[$\mathsection 1.9$]{Br83}. 
Brieskorn~\cite[Tabelle 4, p. 68]{Br83} gives a data that induces their discriminant quadratic forms for  hyperbolic lattices of type $L''_{p,q,r}$. 
Following this result, we reproduce the quadratic form of lattices $L''_{p,q,r}$ in Table~\ref{L''_{p,q,r}}. 

\renewcommand{\arraystretch}{1.8}

\begin{tiny}
\[
\begin{array}{llll l@{\hspace{10mm}} llll}
S & \rho(S) & \delta(S) & q(S) & & S & \rho(S) & \delta(S) & q(S) \\
\cline{1-4}\cline{6-9}
\noalign{\vspace{.5mm}}
\cline{1-4}\cline{6-9}
L''_{2,3,7} & 10 & -1 & 1 & & L''_{2,3,19} & 22 & -13 & w^1_{13,1} \\
\cline{1-4}\cline{6-9}
L''_{2,3,8} & 11 & 2 & w^{-1}_{2,1} & & L''_{3,4,5} & 10 & -13 & w^{-1}_{13,1} \\
\cline{1-4}\cline{6-9}
L''_{2,4,5} & 9 & -2 & w^1_{2,1} & & L''_{2,3,20} & 23 & 14 & w^1_{2,1}\perp w^1_{7,1} \\
\cline{1-4}\cline{6-9}
L''_{2,3,9} & 12 & -3 & w^1_{3,1} & & L''_{2,4,11} & 15 & 14 & w^1_{2,1}\perp w^1_{7,1} \\
\cline{1-4}\cline{6-9}
L''_{3,3,4} & 8 & -3 & w^{-1}_{3,1} & & L''_{2,5,8} & 13 & 14 &  w^{-1}_{2,1}\perp w^{-1}_{7,1} \\ 
\cline{1-4}\cline{6-9}
L''_{2,3,10} & 13 & 4 & w^5_{2,2} & & L''_{2,3,21} & 24 & -15  & w^{-1}_{3,1}\perp w^{-1}_{5,1} \\
\cline{1-4}\cline{6-9}
L''_{2,4,6} & 10 & -4  & w^1_{2,1}\perp w^{-1}_{2,1} & & L''_{3,3,8} & 12 & -15 & w^{-1}_{3,1}\perp w^1_{5,1} \\
\cline{1-4}\cline{6-9}
L''_{2,3,11} & 14 & -5 & w^1_{5,1} & & L''_{2,3,22} & 25 & 16 & w^1_{2,4} \\
\cline{1-4}\cline{6-9}
L''_{2,5,5} & 10 & -5 & w^{-1}_{5,1} & & L''_{2,4,12} & 16 & -16 & w^1_{2,1}\perp w^1_{2,3} \\
\cline{1-4}\cline{6-9}
L''_{2,3,12} & 15 & 6 & w^1_{2,1}\perp w^{-1}_{3,1} & & L''_{2,6,7} & 13 & 16 &  w^5_{2,4} \\
\cline{1-4}\cline{6-9}
L''_{2,4,7} & 11 & 6 & w^1_{2,1}\perp w^1_{3,1} & & L''_{4,4,4} & 10 & -16 & v_2 \\
\cline{1-4}\cline{6-9}
L''_{3,3,5} & 9 & 6 & w^{-1}_{2,1}\perp w^{-1}_{3,1} & & L''_{2,3,23} & 26 & -17 & w^1_{17,1} \\ 
\cline{1-4}\cline{6-9}
L''_{2,3,13} & 16 & -7 & w^1_{7,1} & & L''_{2,5,9} & 14 & -17 & w^{-1}_{17,1} \\ 
\cline{1-4}\cline{6-9}
L''_{2,3,14} & 17 & 8 & w^1_{2,3} & & L''_{2,3,24} & 27 & 18 & w^{-1}_{2,1}\perp w^{-1}_{3,2} \\ 
\cline{1-4}\cline{6-9}
L''_{2,4,8} & 12 & -8 & w^1_{2,1}\perp w^5_{2,2} & & L''_{2,4,13} & 17 & 18 & w^1_{2,1}\perp w^1_{3,2} \\
\cline{1-4}\cline{6-9}
L''_{2,5,6} & 11 & 8 & w^{-5}_{2,3} & & L''_{3,3,9} & 13 & 18 &  w^1_{2,1}\perp (w^1_{3,1})^2 \\
\cline{1-4}\cline{6-9}
L''_{3,4,4} & 9 & 8  & w^5_{2,3} & & L''_{3,4,6} & 11 & 18 &  w^{-1}_{2,1}\perp w^1_{3,2} \\ 
\cline{1-4}\cline{6-9}
L''_{2,3,15} & 18 & -9 & w^1_{3,2} & & L''_{2,3,25} & 28 & -19 & w^1_{19,1} \\
\cline{1-4}\cline{6-9}
L''_{3,3,6} & 10 & -9 & w^{-1}_{3,1}\perp w^1_{3,1} & & L''_{2,3,26} & 29 & -20 & w^1_{2,2}\perp w^1_{5,1} \\
\cline{1-4}\cline{6-9}
L''_{2,3,16} & 19 & 10 & w^{-1}_{2,1}\perp w^{-1}_{5,1} & &  L''_{2,4,14} & 18 & -20 & w^1_{2,1}\perp w^{-1}_{2,1}\perp w^{-1}_{5,1} \\
\cline{1-4}\cline{6-9}
L''_{2,4,9} & 13 & 10 &  w^1_{2,1}\perp w^1_{5,1} & & L''_{2,5,10} & 15 & 20 &  w^{-5}_{2,2}\perp w^{-1}_{5,1} \\ 
\cline{1-4}\cline{6-9}
L''_{2,3,17} & 20 & -11 & w^1_{11,1} & & L''_{2,6,8} & 14 & -20 &  w^1_{2,1}\perp w^{-1}_{2,1}\perp w^1_{5,1} \\
\cline{1-4}\cline{6-9}
L''_{2,5,7} & 12 & -11 & w^1_{11,1} & & L''_{3,5,5} & 11 & 20 &  w^{-1}_{2,2}\perp w^{-1}_{5,1} \\
\cline{1-4}\cline{6-9}
L''_{2,3,18} & 21 & 12 & w^{-1}_{2,2}\perp w^1_{3,1} & & L''_{2,3,27} & 30 & -21 & w^1_{3,1}\perp w^{-1}_{7,1} \\
\cline{1-4}\cline{6-9}
L''_{2,4,10} & 14 & -12 & (w^1_{2,1})^2\perp w^{-1}_{3,1} & & L''_{2,7,7} & 14 & -21 & w^{-1}_{3,1}\perp w^1_{7,1} \\
\cline{1-4}\cline{6-9}
L''_{2,6,6} & 12 & -12  & u_1\perp w^1_{3,1} & & L''_{3,3,10} & 14 & -21 & w^{-1}_{3,1}\perp w^1_{7,1} \\
\cline{1-4}\cline{6-9}
L''_{3,3,7} & 11 & 12 & w^5_{2,2}\perp w^{-1}_{3,1} \\
\cline{1-4}
\end{array}
\]
\begingroup
\captionof{table}{Hyperbolic lattices of type $L''_{p,q,r}$}\label{L''_{p,q,r}}
\endgroup 
\end{tiny}

\bigskip
Let $(\rho,\,\delta)$ be a pair that appears as the rank and discriminant number of each Picard lattice of the family of weighted $K3$ surfaces associated to simple $K3$ singularities, given by Belcastro~\cite{Be02}. 
We study the set $\mathcal{L}_{(\rho,\,\delta)}$ of non-degenerate even hyperbolic lattices that are primitive sublattices of the $K3$ lattice $\Lambda_{K3}$, of signature $(1,\,\rho-1)$, and discriminant number $\delta$. 
First recall results by Nikulin~\cite{Ni80}. 

\begin{cor}[{\cite[Corollary 1.10.2]{Ni80}}]\label{Existence}
There exists an even lattice of signature $(t_{(+)},\,t_{(-)})$ and with discriminant quadratic form $q$ if the following conditions are simultaneously satisfied: 
\begin{enumerate}
\item[(i)] $t_{(+)}-t_{(-)}\equiv\sign{(q)}\mod8$, \quad 
(ii)\, $t_{(+)}\geq0,\, t_{(-)}\geq0$ and \, $t_{(+)}+t_{(-)}\geq \mathfrak{l}(G_q)$. $\square$
\end{enumerate}
\end{cor}

\begin{cor}[{\cite[Corollary 1.13.3]{Ni80}}]\label{Uniqueness}
There exists a unique even lattice of signature $(t_{(+)},\,t_{(-)})$ and with discriminant quadratic form $q$ if the following three conditions are simultaneously satisfied: 
\begin{enumerate}
\item[(i)] $t_{(+)}-t_{(-)}\equiv\sign{(q)}\mod8$, \quad
(ii) $t_{(+)}\geq1,\, t_{(-)}\geq1$, \quad (iii) $t_{(+)}+t_{(-)}\geq 2+ \mathfrak{l}(G_q)$. $\square$
\end{enumerate}
\end{cor}

\subsection{Main results}
We describe the set $\mathcal{L}_{(\rho,\,\delta)}$ in special cases. 
\begin{thm}\label{LatticeSimple}
There exists a unique lattice with signature $(1,\,\rho-1)$ and of discriminant number $\delta$ up to isomorphism, if a pair $(\rho,\,\delta)$ is given with 
\begin{itemize}
\item[$(1)$] $\delta$ is a prime number, as in Table~\ref{Prime}, 
\item[$(2)$] $\rho=1$, as in Table~\ref{rho1}, \quad 
$(3)$ $\delta=1$, as in Table~\ref{delta1}. 
\end{itemize}
\end{thm}
\noindent
\proof
For all pairs $(\rho,\,\delta)$, one can find all quadratic forms $q$ with signature $2-\rho$ and discriminant groups $G$ of order $\delta$. 
With the data, one can compare $\rho$ and the values $l(G)$, and $2+l(G)$, respectively. \\
(1):\, If $\rho=2$, then, the statement follows from an explicit computation to give the lattice. 
If $\rho\not=2$, then, Corollaries \ref{Existence} and \ref{Uniqueness} can apply and thus the statement is verified. \\
(2) and (3):\, It is clear that Corollaries \ref{Existence} and \ref{Uniqueness} can apply and thus the statement is verified. 
$\square$

In Tables~\ref{rho1},~\ref{delta1} and~\ref{Prime}, we present the ranks $\rho$ and discriminant numbers $\delta$ as is given in ~\cite{Be02}, then, the discriminant group $G$ and quadratic form $q$ on $G$, and lattices. 
In the sixth column, we distinguish the family in Yonemura's classification. 

\renewcommand{\arraystretch}{1.8}

\begin{tiny}
\begin{tabular}{ccc}
\begin{minipage}{7cm}
\[
\begin{array}{lrllll}
\rho & \delta & G & q & \textnormal{\tiny{Lattice}} & \textnormal{\tiny{No.} }
\\ \hline \noalign{\vspace{.3mm}}\hline
1 & 4 & \mathbb{Z}_4 & w^1_{2,2} & \langle 4\rangle & 1  \\
\hline
1 & 2 & \mathbb{Z}_2 & w^1_{2,1} & \langle 2\rangle & 5  \\
\hline
\end{array} 
\]
\begingroup
\captionof{table}{Lattices of rank one}\label{rho1} 
\endgroup
\end{minipage}
& & 
\begin{minipage}{7cm}
\[
\begin{array}{lrlllr}
\rho & \delta & G & q & \textnormal{\tiny{Lattice}} & \textnormal{\tiny{No.} }
\\ \hline \noalign{\vspace{.3mm}}\hline
 2  & 1 & \{0\} & 1 & U & 10   \\
\hline
10  & 1 & \{ 0\} & 1 & U\oplus E_8 & 14, 28, 45, 51  \\
\hline
18  & 1 & \{ 0\} & 1 & U\oplus E_8^{\oplus2} & 46, 65, 80 \\
\hline
\end{array}
\]
\begingroup
\captionof{table}{Unimodular lattices}\label{delta1}
\endgroup 
\end{minipage}
\end{tabular}

\[
\begin{array}{lrlllr}
\rho & \delta & G & q & \textnormal{Lattice} & \textnormal{No.}
\\ \hline \noalign{\vspace{.3mm}}\hline
2  &  5 & \mathbb{Z}_5 & w^{-1}_{5,1} & \left( \mathbb{Z}^2,\, \begin{pmatrix} 2 & 1 \\ 1 & -2 \end{pmatrix}\right) & 21  \\
\hline
3  &  2 & \mathbb{Z}_2 & w^{-1}_{2,1} & U\oplus A_1 & 42  \\
\hline
4 &  3 & \mathbb{Z}_3 & w^1_{3,1} & U\oplus A_2 & 25  \\
\cline{2-6}
   &  7 & \mathbb{Z}_7 & w^{-1}_{7,1} & U\oplus S_{w^{-1}_{7,1}} & 66  \\
\hline
8  &  3 & \mathbb{Z}_3 & w^{-1}_{3,1} & U\oplus E_6 & 13, 72  \\
\cline{2-6}
  &  11 & \mathbb{Z}_{11} & w^{-1}_{11,1} & U\oplus S_{w^{-1}_{11,1}} & 89 \\
\hline
9  &  2 & \mathbb{Z}_2 & w^1_{2,1} & U\oplus E_7 & 50, 82 \\
\hline
10 & 5 & \mathbb{Z}_5 & w^{-1}_{5,1} & L''_{2,5,5} & 9, 71 \\
\cline{2-6}
    & 13 & \mathbb{Z}_{13} & w^{-1}_{13,1} & L''_{3,4,5} & 87 \\
\hline
11 &  2 & \mathbb{Z}_2 & w^{-1}_{2,1} & U\oplus A_1\oplus E_8 & 38, 77  \\
\hline
12 &  3 & \mathbb{Z}_3 & w^1_{3,1} & U\oplus A_2\oplus E_8 & 20, 59  \\
\hline
14 &  17 & \mathbb{Z}_{17} & w^{-1}_{17,1} & L''_{2,5,9} & 95 \\ 
\hline
16 &  3 & \mathbb{Z}_3 & w^{-1}_{3,1} & U\oplus E_6\oplus E_8 & 43, 48, 88   \\
\cline{2-6}
 &  7 & \mathbb{Z}_7 & w^1_{7,1} & U\oplus A_6\oplus E_8 & 35 \\
\cline{2-6}
 &  19 & \mathbb{Z}_{19} & w^{-1}_{19,1} & U\oplus S_{w^{-1}_{19,1}} & 94 \\
\hline
17 &  2 & \mathbb{Z}_2 & w^1_{2,1} & U\oplus E_7\oplus E_8 & 68, 83, 92  \\
\hline
18 &  5 & \mathbb{Z}_5 & w^{-1}_{5,1} & L''_{2,5,5}\oplus E_8 & 30, 86   \\
\hline
19 &  2 & \mathbb{Z}_2 & w^{-1}_{2,1} & U\oplus A_1\oplus E_8^{\oplus2} & 56, 73   \\
\hline
\end{array}
\]
\begingroup
\captionof{table}{Lattices with prime discriminant numbers}\label{Prime}
\endgroup
\end{tiny}

\bigskip


We deal with $\mathcal{L}_{(\rho,\,\delta)}$ for the remaining cases. 

\begin{thm}\label{MainThmLattice}
For a given pair $(\rho,\,\delta)$ as in Table~\ref{Reclassify}, one can give all quadratic forms of signature $2-\rho$ and defined on the discriminant group of order $\delta$ up to isomorphism. 
Moreover, for each triple $(\rho,\,\delta,\, q)$, there exists a unique lattice of signature $(1,\,\rho-1)$ that admits $q$ as its discriminant quadratic form. 
\end{thm}

In Table~\ref{Reclassify}, one finds in the first and second columns $\rho$ and $\delta$ that we treat, in the third and fourth columns are the discriminant group $G$ and the quadratic form $q$, in the fifth column, the lattices, and in the sixth column, the index of the corresponding weight system.  
Note that the signature of a quadratic form $q$ is given by $\sign{(q)}=2-\rho\mod8$. 
\\

\noindent
{\sc Proof of Theorem~\ref{MainThmLattice} }
We take three steps. 

{\it Step 1. } We prove the following lemma. 
\begin{lem}\label{MainLem}
For a pair $(\rho,\,\delta)$, there exist discriminant quadratic forms of signature $2{-}\rho$ that are defined on the discriminant group of order $\delta$. 
The results are given in Table~\ref{Reclassify}, columns 3 and 4. 
\end{lem}
\noindent
{\sc Proof of Lemma~\ref{MainLem} }
One can get all quadratic forms of signature ${2-\rho}$ and the discriminant groups by using our previous study. 

We give an example $(\rho,\,\delta)=(3,4)$ for what we excluded from a set of candidates of quadratic forms. 
One may expect to take a quadratic form which is defined on a finite group $\mathbb{Z}_2^2$. 
However, it is in fact impossible:  
indeed, quadratic forms on the group $\mathbb{Z}_2^2$ are 
\[
(w^1_{2,1})^2,\, w^1_{2,1}\perp w^{-1}_{2,1},\, (w^{-1}_{2,1})^2,\, u_1,\quad \textnormal{or} \quad v_1, 
\]
and the signatures of the forms are respectively
\[
2, \, 0, \, 6, \, 0, \quad \textnormal{or} \quad 4. 
\]
So, any lattice of signature $(1,2)$ does not admit a quadratic form of signature equal to $2-\rho=2-3\equiv7\mod8$. 
QED for Lemma~\ref{MainLem}.

{\it Step 2. } 
By {\it Step 1}, it is verified that they satisfy the condition (1) in Corollaries \ref{Existence} and \ref{Uniqueness}, and can compare the values $\rho$ and $\mathfrak{l}(G)$, $2+\mathfrak{l}(G)$ explicitely. 
Thus, one can verify that they also satisfy the condition (2) in Corollaries \ref{Existence} and \ref{Uniqueness}. 
Therefore, by Corollaries \ref{Existence} and \ref{Uniqueness}, we can conclude that a lattice with invariant $(\rho,\,\delta,\, q)$ uniquely exists. 

{\it Step 3. }
By a case-by-case argument, one can determine a lattice of signature $(1,\,\rho-1)$ that admits $q$ as its discriminant quadratic form. 
By studying~\cite{Be02}, one can distinguish the weight system with obtained lattices. 
$\square$ 

\begin{remark}
From Table~\ref{Reclassify}, we have confirmed that only a pair of rank and discriminant does not necessarily determine a lattice, and have observed that in case of families of weighted $K3$ surfaces, the Picard lattice of the family is determined (as in~\cite{Be02}) due to an embedding of surfaces into the weighted projective space. 
\end{remark}

Finally, we consider $K3$ surfaces $S=S_L$ polarized by a lattice $L\in\mathcal{L}_{(\rho,\,\delta)}$ such that $\Pic{(S)} = L$. 
Denote by $\mathcal{M}_{\rho,\,\delta}$ the collection of such surfaces. 

\begin{prop}
For any $(\rho,\,\delta)$ in Tables~\ref{rho1},~\ref{delta1}~\ref{Prime}, and~\ref{Reclassify}, all lattices in $\mathcal{L}_{(\rho,\,\delta)}$ represent zero except $\rho=1$ and $2$.  
\end{prop}
\noindent
\proof
If $\rho\geq5$, all lattices in $\mathcal{L}_{(\rho,\,\delta)}$ represent zero~\cite{SeLec}. 
If $\rho<5$, the statement is verified by the tables, except two cases:
$\left( \mathbb{Z}^2,\, \begin{pmatrix} 2 & 1 \\ 1 & -2\end{pmatrix}\right)$ in $\mathcal{L}_{(2,5)}$, 
and 
$\langle 2\rangle\oplus A_1\oplus A_2 $ in $\mathcal{L}_{(4,12)}$. 

Let $\langle g_1,\, g_2\rangle_{\mathbb{Z}}$ be a basis for $\left( \mathbb{Z}^2,\, \begin{pmatrix} 2 & 1 \\ 1 & -2\end{pmatrix}\right)$ such that $g_1^2=2$, $g_2^2=-2$, and $g_1.g_2 = 1$. 

An element $a g_1 + b g_2$ in the lattice ($a,\,b\in\mathbb{Z}\slash5\mathbb{Z}$) is of norm zero if and only if $a=b=0$ since $(a g_1 + b g_2)^2 = 2(a^2+ab-b^2)$. 
Therefore, the lattice does not represent zero. 

Let $\langle g_1,\, g_2,\, g_3,\, g_4\rangle_{\mathbb{Z}}$ be a basis for $\langle 2\rangle\oplus A_1\oplus A_2 $ such that $g_1^2=2$, $g_2^2=-2$, and $\langle g_3,\, g_4\rangle_\mathbb{Z} = A_2$. 
Then, we have $(g_1 + g_2)^2 = 0$ thus the lattice represents zero. \QED

\bigskip

Since the class of divisor of self-intersection zero can be represented by an elliptic curve, we can conclude: 
\begin{cor}\label{AllElliptic}
Let $\rho\not=1,\,2$. 
All $K3$ surfaces in $\mathcal{M}_{(\rho,\,\delta)}$ polarized by lattices in $\mathcal{L}_{(\rho,\,\delta)}$ have a structure of elliptic surface. \QED
\end{cor}

\bigskip

\renewcommand{\arraystretch}{1.8}
\begin{tiny}
\[
\begin{array}{crllll}
\rho & \delta & G & q & \textnormal{Lattice} & \textnormal{No. }
\\ \hline \noalign{\vspace{.3mm}}\hline
3 & 4 & \mathbb{Z}_4 & w^{-1}_{2,2} & U\oplus\langle-4\rangle &7  \\
\cline{3-6}
\hline
4 & 12 & \mathbb{Z}_2\times\mathbb{Z}_6 & w^1_{2,1}\perp w^{-1}_{2,1}\perp w^1_{3,1} & \langle 2\rangle\oplus A_1\oplus A_2 \\
\cline{4-6}
 &  & & v_1\perp w^{-1}_{3,1} & U\oplus S_{v_1\perp w^{-1}_{3,1}} & 3 \\
\hline
6 & 4 & \mathbb{Z}_2^2 & v_1 & U\oplus D_4 & 12 \\
\cline{2-6}
 & 16 & \mathbb{Z}_2^4 & u_1\perp v_1 & U(2)\oplus D_4 & 6 \\
\cline{3-6}
 &     & \mathbb{Z}_4^2 & w^1_{2,2}\perp w^{-5}_{2,2} & \langle4\rangle\oplus D_5\\ 
\cline{4-6}
 &     & & w^{-1}_{2,2}\perp w^5_{2,2} & U\oplus\langle-4\rangle\oplus A_3 \\
\hline
7 & 4 & \mathbb{Z}_4 & w^{-5}_{2,2} & U\oplus D_5 & 44 \\
\cline{2-6}
 & 8 & \mathbb{Z}_2^3 & v_1\perp w^{-1}_{2,1} & U\oplus A_1\oplus D_4 & 40 \\
\cline{2-6}
 & 12 & \mathbb{Z}_3\times\mathbb{Z}_4 & w^1_{3,1}\perp w^5_{2,2} & U\oplus A_2\oplus A_3 & \\
\cline{4-6}
 & & & w^{-1}_{3,1}\perp w^1_{2,2} & \langle4\rangle\oplus E_6 & 8  \\
\cline{2-6}
 &  32   &  \mathbb{Z}_2^2\times\mathbb{Z}_8 & v_1\perp w^{-1}_{2,3} & U\oplus\langle-8\rangle\oplus D_4 & 19 \\ 
\cline{4-6}
 & &  & v_1\perp w^{-5}_{2,3} & D_4\oplus S_{w^{-5}_{2,3}} & \\ 
\cline{4-6}
 & & & (w^1_{2,1})^2\perp w^1_{2,3} & \langle8\rangle\oplus D_6\\
\cline{3-6}
 &     & \mathbb{Z}_2^5 & (w^{-1}_{2,1})^5 & U\oplus A_1^{\oplus 5} \\ 
\cline{3-6}
 &     & \mathbb{Z}_2\times\mathbb{Z}_4^2 & w^1_{2,1}\perp (w^1_{2,2})^2 & \langle2\rangle\oplus A_3^{\oplus2}\\
\cline{4-6}
 &     & & w^1_{2,1}\perp w^{-1}_{2,2}\perp w^{-5}_{2,2} & \langle2\rangle\oplus\langle-4\rangle\oplus D_5\\
\cline{4-6}
 &     & & w^{-1}_{2,1}\perp w^{-1}_{2,2}\perp w^{5}_{2,2} & U\oplus\langle-4\rangle\oplus A_1\oplus A_3\\
\hline
8 & 12 & \mathbb{Z}_3\times\mathbb{Z}_4 & v_1\perp w^1_{3,1} & U\oplus A_2\oplus D_4 & 24 \\
\cline{3-6}
 &  & \mathbb{Z}_2\times\mathbb{Z}_6 & w^{-1}_{2,1}\perp w^1_{2,1}\perp w^{-1}_{3,1} & \langle2\rangle\oplus A_1\oplus E_6 \\
\cline{2-6}
 & 20 & \mathbb{Z}_2\times\mathbb{Z}_{10} & (w^1_{2,1})^2\perp w^{-1}_{5,1} & S_{w^{-1}_{5,1}}\oplus D_6 & 63 \\
\cline{4-6}
 & & & (w^{-1}_{2,1})^2\perp w^1_{5,1} & U\oplus A_1^{\oplus2}\oplus A_4\\
\cline{2-6}
 & 27 & \mathbb{Z}_3\times\mathbb{Z}_9 & w^{-1}_{3,1}\perp w^{-1}_{3,2} & E_6\oplus S_{w^{-1}_{3,2}} & 18 \\
\cline{4-6}
 &     & & w^{-1}_{3,1}\perp w^1_{3,2} & S_{w^{1}_{3,2}}\oplus E_6\\
\cline{3-6}
 &     & \mathbb{Z}_3^3 & (w^1_{3,1})^3 & U\oplus A_2^{\oplus 3} \\
\cline{3-6}
 &     & \mathbb{Z}_{27} & w^{-1}_{3,3} & U\oplus S_{w^{-1}_{3,3}} \\
\hline
9 & 6 & \mathbb{Z}_2\times\mathbb{Z}_3 & w^{-1}_{2,1}\perp w^{-1}_{3,1} & U\oplus A_1\oplus E_6 & 39 \\
\cline{2-6}
 & 8 & \mathbb{Z}_8 & w^1_{2,3} & U\oplus A_7 & \\
\cline{4-6}
 &    &            & w^5_{2,3} & L''_{3,4,4} & 37 \\
\cline{3-6}
 &    & \mathbb{Z}_2^3  & (w^1_{2,1})^2\perp w^{-1}_{2,1} & \langle2\rangle\oplus A_1\oplus E_7\\
\hline
\end{array}
\]
\hfill (Table~\ref{Reclassify}. Continued on the next page. )

\[
\begin{array}{crlllr}
\hline
10 & 4 & \mathbb{Z}_2^2 & w^1_{2,1}\perp w^{-1}_{2,1} & U\oplus A_1\oplus E_7 & 78 \\
\cline{4-6}
 & & & u_1 & U\oplus D_8 & \\
\cline{2-6}
 & 9 & \mathbb{Z}_9 & w^1_{3,2} & U\oplus A_8 & \\
\cline{4-6}
 & & & w^{-1}_{3,2} & S_{w^{-1}_{3,2}}\oplus E_8 \\
\cline{3-6}
 & & \mathbb{Z}_3^2 & w^1_{3,1}\perp w^{-1}_{3,1} & U\oplus A_2\oplus E_6 & 22 \\
\cline{2-6}
 & 16 & \mathbb{Z}_4^2 & v_2 & L''_{4,4,4} & 4 \\
\cline{4-6}
 & & & w^5_{2,2}\perp w^{-5}_{2,2} & U\oplus A_3\oplus D_5 & \\ 
\cline{4-6}
 & & & w^1_{2,2}\perp w^{-1}_{2,2} & \langle4\rangle\oplus\langle-4\rangle\oplus E_8\\
\cline{3-6}
\cline{3-6}
 &  & \mathbb{Z}_2^4 & (v_1)^2 & U\oplus D_4^{\oplus2} & \\ 
\cline{4-6}
 &  & & u_1\perp w^1_{2,1}\perp w^{-1}_{2,1} & \langle2\rangle\oplus A_1\oplus D_8\\
\cline{3-6}
 &     & \mathbb{Z}_2\times\mathbb{Z}_8 & w^{-1}_{2,1}\perp w^1_{2,3} & U\oplus A_1\oplus A_7 & \\
\cline{4-6}
\cline{4-6}
 &  & & w^1_{2,1}\perp w^{-1}_{2,3} & \langle2\rangle\oplus\langle-8\rangle\oplus E_8 \\ 
\cline{2-6}
 & 64 & \mathbb{Z}_2^6 & (u_1)^3 & U(2)\oplus D_4^{\oplus 2} & 32 \\
\cline{3-6}
 &     & \mathbb{Z}_2^3\times\mathbb{Z}_8 & w^1_{2,1}\perp (w^{-1}_{2,1})^2\perp w^1_{2,3} & \langle2\rangle\oplus A_1^{\oplus2}\oplus A_7 \\ 
\cline{4-6}
 &     & & (w^1_{2,1})^2\perp w^{-1}_{2,1}\perp w^{-1}_{2,3} & U\oplus\langle-8\rangle\oplus A_1\oplus D_6\\
\cline{3-6}
& & \mathbb{Z}_2^2\times\mathbb{Z}_4^2  & (w^{-1}_{2,1})^2\perp (w^5_{2,2})^2 & U\oplus A_1^{\oplus 2}\oplus A_3^{\oplus2}   \\
\cline{4-6}
 &     &  & v_1\perp w^{-1}_{2,2}\perp w^{5}_{2,2} & U\oplus\langle-4\rangle\oplus A_3\oplus D_4\\
\cline{3-6}
 &     & \mathbb{Z}_2\times\mathbb{Z}_{32} & w^{1}_{2,1}\perp w^{-1}_{2,5} & U\oplus\langle-32\rangle\oplus E_7 \\
\cline{4-6}
 &     & & w^{1}_{2,1}\perp w^{-5}_{2,5} & S_{w^{-5}_{2,5}}\oplus E_7\\
\cline{4-6}
 &     & & w^{-1}_{2,1}\perp w^{1}_{2,5} & U\oplus A_1\oplus S_{w^{1}_{2,5}} \\
 \cline{4-6}
 &     & & w^{-1}_{2,1}\perp w^{5}_{2,5} & U\oplus A_1\oplus S_{w^{5}_{2,5}} \\
\cline{3-6}
 &     & \mathbb{Z}_4\times\mathbb{Z}_{16} & w^{-1}_{2,2}\perp w^{1}_{2,4} & \langle16\rangle\oplus D_9\\
\cline{4-6}
 &     &  & w^{1}_{2,2}\perp w^{-1}_{2,4} & U\oplus\langle-16\rangle\oplus D_7\\
\cline{3-6}
 &     & \mathbb{Z}_8^2 & w^{1}_{2,3}\perp w^{-5}_{2,3} & S_{w^{-5}_{2,3}}\oplus A_7 \\
\cline{4-6}
 & & & w^{1}_{2,3}\perp w^{-1}_{2,3} & U\oplus\langle-8\rangle\oplus A_7\\
\hline
\end{array}
\]
\hfill (Table~\ref{Reclassify}. Continued on the next page. )

\[
\begin{array}{crlllr}
\hline
11 & 6 & \mathbb{Z}_6 & w^1_{2,1}\perp w^1_{3,1} & U\oplus E_7\oplus A_2 & 60 \\
\cline{2-6}
 & 8 &  \mathbb{Z}_8 & w^{-5}_{2,3} & L''_{2,5,6} & 58  \\
\cline{3-6}
 &    & \mathbb{Z}_2^3 & w^1_{2,1}\perp (w^{-1}_{2,1})^2 & U\oplus A_1^{\oplus 2}\oplus E_7 & \\
\cline{2-6}
 & 64 & \mathbb{Z}_2^4\times\mathbb{Z}_4 & (v_1)^2\perp w^{-1}_{2,2} & U\oplus\langle-4\rangle\oplus D_4^{\oplus2} & 23 \\
\cline{4-6}
 & & & (w^{-1}_{2,1})^4\perp w^{-5}_{2,2} & U\oplus A_1^{\oplus4}\oplus D_5\\
\cline{4-6}
 & & & v_1\perp(w^{-1}_{2,1})^2\perp w^5_{2,2} & U\oplus A_1^{\oplus2}\oplus A_3\oplus D_4 \\
\cline{4-6}
 & & & (w^{1}_{2,1})^2\perp (w^{-1}_{2,1})^2\perp w^{-1}_{2,2} & \langle2\rangle\oplus\langle-4\rangle\oplus A_1^{\oplus2}\oplus E_7\\
\cline{3-6}
 &     & \mathbb{Z}_4^3 & (w^5_{2,2})^3 & U\oplus A_3^{\oplus3} &  \\
\cline{4-6}
 &     &  & w^{-1}_{2,2}\perp w^{5}_{2,2}\perp w^{-5}_{2,2} & U\langle-4\rangle\oplus A_3\oplus D_5\\
\cline{3-6}
 &     & \mathbb{Z}_2^2\times\mathbb{Z}_{16} & (w^{1}_{2,1})^2\perp w^{5}_{2,4} & U\oplus D_6\oplus S_{w^{5}_{2,4}} \\ 
\cline{4-6}
 &     &  & (w^{-1}_{2,1})^2\perp w^{1}_{2,4} & \langle16\rangle\oplus A_1^{\oplus2}\oplus E_8\\
\cline{4-6}
 &     &  & v_1\perp w^{-5}_{2,4} & U\oplus D_4\oplus S_{w_{2,4}^{-5}}  \\
\cline{4-6}
 &     &  & u_1\perp w^{-1}_{2,4} & U\oplus\langle-16\rangle\oplus D_8\\ 
\cline{3-6}
 &     & \mathbb{Z}_2\times\mathbb{Z}_4\times\mathbb{Z}_8 & w^{-1}_{2,1}\perp w^{-1}_{2,2}\perp w^1_{2,3} & U\oplus \langle-4\rangle\oplus A_1\oplus A_7 \\ 
\cline{3-6}
 &     & \mathbb{Z}_{64} & w^{-1}_{2,6} & U\oplus\langle-64\rangle\oplus E_8\\
\hline
12 & 12 & \mathbb{Z}_{12} & u_1\perp w^1_{3,1} & U\oplus A_2\oplus D_8 &  \\
\cline{4-6}
  &  &  & v_1\perp w^{-1}_{3,1} & U\oplus E_6\oplus D_4 & 11 \\
\cline{3-6}
    &      & \mathbb{Z}_2\times\mathbb{Z}_6 & w^1_{2,1}\perp w^{-1}_{2,1}\perp w^1_{3,1} & U\oplus E_7\oplus A_1\oplus A_2 & \\
\cline{2-6}
 & 48 & \mathbb{Z}_2^4\times\mathbb{Z}_3 & v_1\perp w^1_{2,1}\perp w^{-1}_{2,1}\perp w^{-1}_{3,1} & \langle2\rangle\oplus A_1\oplus D_4\oplus E_6 & 33 \\
\cline{3-6}
 &     & \mathbb{Z}_2^3\times\mathbb{Z}_6 &  (w^{-1}_{2,1})^4\perp w^{-1}_{3,1} & U\oplus A_1^{\oplus4}\oplus E_6 \\
\cline{4-6}
 &     &   & (w^1_{2,1})^2\perp (w^{-1}_{2,1})^2\perp w^1_{3,1} & \langle2\rangle\oplus A_1^{\oplus2}\oplus A_2\oplus E_7 \\
\cline{3-6}
 &     & \mathbb{Z}_2^2\times\mathbb{Z}_{12} & (w^1_{2,1})^2\perp (w^{-1}_{2,1})^2\perp w^1_{3,1}  & U\oplus A_1^{\oplus2}\oplus A_2\oplus D_6 &  \\
\cline{3-6}
 & & \mathbb{Z}_6\times\mathbb{Z}_8 & w^{-1}_{2,1}\perp w^1_{2,3}\perp w^1_{3,1} & U\oplus A_1\oplus A_2\oplus A_7 &  \\
\cline{4-6}
 & &  & w^1_{2,1}\perp w^{-5}_{2,3}\perp w^{1}_{3,1} & A_2\oplus E_7\oplus S_{w^{-5}_{2,3}}\\
\cline{3-6}
 &      & \mathbb{Z}_4\times\mathbb{Z}_{12} & w^{1}_{2,2}\perp w^{-1}_{2,2}\perp w^1_{3,1} & \langle4\rangle\oplus\langle-4\rangle\oplus A_2\oplus E_8\\
\cline{4-6}
 &      &  & w^{1}_{2,2}\perp w^{-5}_{2,2}\perp w^{-1}_{3,1} & \langle4\rangle\oplus D_5\oplus E_6 \\ 
\cline{4-6}
 &      &  & w^{-1}_{2,2}\perp w^{5}_{2,2}\perp w^{-1}_{3,1} & U\oplus\langle-4\rangle\oplus A_3\oplus E_6\\
\cline{4-6}
 &      & & w^5_{2,2}\perp w^{-5}_{2,2}\perp w^1_{3,1} & U\oplus A_2\oplus A_3\oplus D_5 &  \\
\cline{2-6}
 & 108  & \mathbb{Z}_3\times\mathbb{Z}_6^2 & v_1\perp w^1_{3,1}\perp (w^{-1}_{3,1})^2 & U(3)\oplus D_4\oplus E_6 & 2 \\
\cline{4-6}
 &        &  & w^{1}_{2,1}\perp w^{-1}_{2,1}\perp w^1_{3,1}\perp w^{-1}_{3,2} & \begin{pmatrix} -2 & 1 \\ 1 & 4\end{pmatrix}\oplus A_1\oplus A_2\oplus E_7\\
\cline{4-6}
 &        &  & w^{1}_{2,1}\perp w^{-1}_{2,1}\perp w^1_{3,1}\perp w^{1}_{3,2} & \langle2\rangle\oplus A_1\oplus A_2\oplus A_8\\
\cline{4-6}
 &        &  & w^{1}_{2,1}\perp w^{-1}_{2,1}\perp w^1_{3,3} & L''_{2,4,6}\oplus\begin{pmatrix} -2 & 1 \\ 1 &  -14\end{pmatrix}\\
\hline
\end{array}
\]
\hfill (Table~\ref{Reclassify}. Continued on the next page. )

\[
\begin{array}{crlllr}
\hline
13 & 8 & \mathbb{Z}_2^3 & (w^{-1}_{2,1})^3 & U\oplus A_1^{\oplus3}\oplus E_8 & 81 \\
\cline{2-6}
 & 12 & \mathbb{Z}_{12} & w^{-1}_{3,1}\perp w^{-5}_{2,2} & U\oplus E_6\oplus D_5 & 41 \\
\cline{4-6}
 &      & & w^1_{3,1}\perp w^{-1}_{2,2} & U\oplus A_2\oplus D_9 &  \\
\cline{2-6}
 & 20 & \mathbb{Z}_{20} & w^1_{2,2}\perp w^1_{5,1} & U\oplus A_4\oplus D_7 & \\
\cline{4-6}
 & & & w^5_{2,2}\perp w^{-1}_{5,1} & L''_{2,5,5}\oplus A_3 & 36 \\
\cline{2-6}
 & 32 & \mathbb{Z}_8\times\mathbb{Z}_2^2 & w^1_{2,3}\perp v_1 & U\oplus A_7\oplus D_4 & 69 \\
\cline{4-6}
 &      &                            & (w^{-1}_{2,1})^2\perp w^{-1}_{2,3} & U\oplus\langle-8\rangle\oplus D_{10}\\
\cline{4-6}
 &      &                            & v_1\perp w^{5}_{2,3} & L''_{3,4,4}\oplus D_4\\
\cline{3-6}
 &  & \mathbb{Z}_2^5 & w^1_{2,1}\perp (w^{-1}_{2,1})^4 & U\oplus A_1^{\oplus4}\oplus E_7 & 75 \\
\cline{3-6}
 & & \mathbb{Z}_2^2\times\mathbb{Z}_8 & v_1\perp w^1_{2,3} & U\oplus A_7\oplus D_4 & \\ 
\cline{3-6}
 &     & \mathbb{Z}_2\times\mathbb{Z}_4^2 & w^{-1}_{2,1}\perp (w^{-5}_{2,2})^2 & U\oplus A_1\oplus D_5^{\oplus2} & \\
\cline{2-6}
 & 56 & \mathbb{Z}_2^3\times\mathbb{Z}_7 & u_1\perp w^{-1}_{2,1}\perp w^{-1}_{7,1} &  U\oplus S_{w^{-1}_{7,1}}\oplus A_1\oplus D_8 \\
\cline{4-6}
 & & & v_1\perp w^{-1}_{2,1}\perp w^1_{7,1} & U\oplus A_1\oplus A_6\oplus D_4 & 85 \\
\cline{4-6}
 & & & w^{-1}_{2,3}\perp w^{-1}_{7,1} & U\oplus\langle-8\rangle\oplus S_{w^{-1}_{7,1}}\oplus E_8\\
\cline{4-6}
 & & & w^{-5}_{2,3}\perp w^{-1}_{7,1} & L''_{2,5,6}\oplus S_{w^{-1}_{7,1}} \\
\hline
14 & 4 & \mathbb{Z}_2^2 & v_1 & U\oplus E_8\oplus D_4 & 27, 49 \\
\cline{3-6}
\cline{2-6}
 & 9 & \mathbb{Z}_3^2 & (w^{-1}_{3,1})^2 & U\oplus A_2^{\oplus2}\oplus E_8 & 67 \\
\cline{2-6}
 & 12 & \mathbb{Z}_2\times\mathbb{Z}_6 & (w^{-1}_{2,1})^2\perp w^1_{3,1} & U\oplus A_1^{\oplus2}\oplus A_2\oplus E_8 & 70 \\
\cline{4-6}
 & & & (w^1_{2,1})^2\perp w^{-1}_{3,1} & U\oplus D_6\oplus E_6\\
\cline{2-6}
 & 16 & \mathbb{Z}_2^4 & u_1\perp v_1 & U\oplus D_4\oplus D_8 & 26, 34, 76 \\
\cline{4-6}
 & & & (w^{-1}_{2,1})^4 & U\oplus A_1^{\oplus4}\oplus E_8\\ 
\cline{3-6}
 &     & \mathbb{Z}_4^2 & w^1_{2,2}\perp w^{-5}_{2,2} & U\oplus D_5\oplus D_7 & \\
\cline{4-6}
 &     &               & w^{-1}_{2,2}\perp w^5_{2,2} & U\oplus A_3\oplus D_9 & \\
\cline{2-6}
 & 25 & \mathbb{Z}_5^2 & w^1_{5,1}\perp w^{-1}_{5,1} & L''_{2,5,5}\oplus A_4 & 17 \\
\cline{2-6}
 & 81 & \mathbb{Z}_3^4 & (w^1_{3,1})^3\perp w^{-1}_{3,1} & U\oplus A_2^{\oplus 3}\oplus E_6 & 15 \\
\cline{3-6}
 & & \mathbb{Z}_3^2\times\mathbb{Z}_9 & (w^1_{3,1})^2\perp w^1_{3,2} & U\oplus A_2^{\oplus2}\oplus A_8 & \\
\cline{3-6}
 &     & \mathbb{Z}_3\times\mathbb{Z}_{27} & w^{1}_{3,1}\perp w^{1}_{3,3} & U\oplus S_{w^1_{3,3}} \oplus A_2\oplus E_8\\
\cline{4-6}
 &     &  & w^{-1}_{3,1}\perp w^{-1}_{3,3} & U\oplus S_{w^{-1}_{3,3}}\oplus E_6 \\
\hline
\end{array}
\]
\hfill (Table~\ref{Reclassify}. Continued on the next page. )

\[
\begin{array}{crlllr}
\hline
15 & 4 & \mathbb{Z}_4 & w^{-5}_{2,2} & U\oplus D_5\oplus E_8 & 79  \\
\cline{2-6}
 & 6 & \mathbb{Z}_6 & w^1_{2,1}\perp w^{-1}_{3,1} & U\oplus E_6\oplus E_7 & 47 \\
\cline{2-6}
 & 16 & \mathbb{Z}_2^2\times\mathbb{Z}_4 & v_1\perp w^{-1}_{2,2} & U\oplus D_4\oplus D_9 & 55 \\
\cline{4-6}
 &  &   & u_1\perp w^{-5}_{2,2} & U\oplus D_5\oplus D_8 & \\
\cline{4-6}
 &  & & (w^{-1}_{2,1})^2\perp w^5_{2,2} & U\oplus A_1^{\oplus2}\oplus A_3\oplus E_8 & \\
\cline{3-6}
 &  & \mathbb{Z}_{16} & w^{-5}_{2,4} & U\oplus S_{w^{-5}_{2,4}}\oplus E_8 \\
\cline{2-6}
 & 24 & \mathbb{Z}_{24} & w^1_{2,3}\perp w^{-1}_{3,1}  & U\oplus A_7\oplus E_6 & 31 \\
\cline{4-6}
 &  &                                 & w^{5}_{2,3}\perp w^{-1}_{3,1} & L''_{3,4,4}\oplus E_6\\
\cline{3-6}
 &  & \mathbb{Z}_2^2\times\mathbb{Z}_6 & (w^1_{2,1})^2\perp w^{-1}_{2,1}\perp w^{-1}_{3,1} & \langle 2\rangle\oplus A_1\oplus E_6\oplus E_7 & \\ 
\cline{4-6}
 & &  & (w^{-1}_{2,1})^3\perp w^1_{3,1} & U\oplus A_1^{\oplus 3}\oplus A_2\oplus E_8 & \\ 
\cline{2-6}
 & 54 & \mathbb{Z}_6\times\mathbb{Z}_9 & w^1_{2,1}\perp w^{-1}_{3,2}\perp w^{-1}_{3,1} & E_6\oplus E_7\oplus S_{w^{-1}_{3,2}} & 53 \\  
\cline{3-6}
 &     & \mathbb{Z}_3\times\mathbb{Z}_{18} & w^1_{2,1}\perp w^{-1}_{3,1}\perp w^1_{3,2} & \langle 2\rangle\oplus A_8\oplus E_6 & \\
\cline{3-6}
 & & \mathbb{Z}_3^2\times\mathbb{Z}_6 & w^1_{2,1}\perp (w^1_{3,1})^3 & U\oplus A_2^{\oplus 3}\oplus E_7 & \\
\cline{3-6}
 & & \mathbb{Z}_2\times\mathbb{Z}_{27} & w^1_{2,1}\perp w^{-1}_{3,3} & U\oplus S_{w^{-1}_{3,3}}\oplus E_7 \\
\hline
16 & 4 & \mathbb{Z}_2^2 & (w^1_{2,1})^2 & U\oplus E_7^{\oplus 2} & 93 \\
\cline{2-6}
 & 16  & \mathbb{Z}_4^2 & w^{-1}_{2,2}\perp w^{-5}_{2,2} & U\oplus D_5\oplus D_9 & 62 \\
\cline{4-6}
 & & & (w^1_{2,2})^2 & U\oplus A_3^{\oplus 2}\oplus E_8 & \\
\cline{3-6}
 & & \mathbb{Z}_2\times\mathbb{Z}_8 & w^1_{2,1}\perp w^1_{2,3} & U\oplus E_7\oplus A_7 & \\
\cline{3-6}
 & & \mathbb{Z}_2^4 & (w^{1}_{2,1})^3\perp w^{-1}_{2,1} & \langle2\rangle\oplus A_1\oplus E_7^{\oplus2}\\
\cline{2-6} 
 & 20 & \mathbb{Z}_2^2\times\mathbb{Z}_5 & (w^1_{2,1})^2\perp w^{-1}_{5,1} & L''_{2,5,5}\oplus D_6 & 29 \\
\cline{4-6}
 & &  & (w^{-1}_{2,1})^2\perp w^1_{5,1} & U\oplus A_1^{\oplus 2}\oplus A_4\oplus E_8 & \\
\cline{2-6}
 & 27 & \mathbb{Z}_3^3 & (w^1_{3,1})^3 & U\oplus A_2^{\oplus 3}\oplus E_8 & 16, 54 \\
\cline{3-6}
 &     & \mathbb{Z}_3\times\mathbb{Z}_9 & w^{-1}_{3,1}\perp w^1_{3,2} & U\oplus A_8\oplus E_6 & \\
\cline{4-6}
 &     &                       & w^{-1}_{3,1}\perp w^{-1}_{3,2} & S_{w^{-1}_{3,2}}\oplus E_6\oplus E_8\\
\cline{3-6}
 &     & \mathbb{Z}_{27} & w^{-1}_{3,3} & U\oplus S_{w^{-1}_{3,3}} \oplus E_8 \\
\hline
17 & 4 & \mathbb{Z}_4 & w^1_{2,2} & U\oplus D_7\oplus E_8 & 64 \\
\cline{2-6}
 & 8 & \mathbb{Z}_2^3 & w^{-1}_{2,1}\perp (w^1_{2,1})^2 & U\oplus A_1\oplus D_6\oplus E_8 & 90 \\
\cline{3-6}
 & & \mathbb{Z}_8 & w^1_{2,3} & U\oplus A_7\oplus E_8 & \\
\cline{4-6} 
  & & & w^5_{2,3} & L''_{3,4,4}\oplus E_8 & 74  \\
\cline{2-6}
 & 12 & \mathbb{Z}_{12} & w^{-5}_{2,2}\perp w^1_{3,1} & U\oplus A_2\oplus D_5\oplus E_8 & 57 \\
\cline{4-6}
 & & & w^{-1}_{2,2}\perp w^{-1}_{3,1} & U\oplus D_9\oplus E_6 & \\
\hline
18 & 4 & \mathbb{Z}_2^2 & w^1_{2,1}\perp w^{-1}_{2,1} & U\oplus E_8\oplus E_7\oplus A_1 & 91 \\
\cline{4-6}
    &    &    & u_1 & U\oplus E_8\oplus D_8 & 61 \\
\cline{2-6}
 & 9 & \mathbb{Z}_9 & w^1_{3,2} & U\oplus E_8\oplus A_8 & 84 \\
\cline{3-6}
 & & \mathbb{Z}_3^2 & w^1_{3,1}\perp w^{-1}_{3,1} & U\oplus E_8\oplus A_2\oplus E_6 & \\
\cline{4-6}
 & & &  w^{-1}_{3,2} & S_{w^{-1}_{3,2}}\oplus E_8^{\oplus2}\\
\hline
19 & 4 & \mathbb{Z}_4 & w^{-1}_{2,2} & U\oplus D_9\oplus E_8 & 52 \\
\hline
\end{array}
\]
\begingroup
\captionof{table}{Lattices with given rank and discriminant number} \label{Reclassify}
\endgroup 
\end{tiny}

\section{A comparison of the Seifert form and Picard lattice}\label{LatticesSingularities}

\subsection{The definition of the Seifert form}
Let a hypersurface singularity $X$ be defined by a polynomial $f$ in variables $(z_0,\, z_1,\ldots,z_n)$. 
On a fibre $F_t:=\varphi^{-1}(\exp{(2\pi it)})$ of the Milnor fibration $\varphi$, there is a homeomorphism $h_t:F=F_0\to F_t$ (see~\cite{Dimca} for more details). 

\begin{defn}
Denote by $lk(\cdot,\cdot)$ the linking number.
Define the {\rm Seifert form} by
\[
L: H_n(F) \times H_n(F)  \rightarrow \mathbb{Z} ; \quad (a,\,b) \mapsto lk(a,\,(h_{1/2})_*b). \qquad\blacksquare
\]
\end{defn}

For the distinguished basis $\mathfrak{D}=(\Delta_i)_{i=1}^\mu$ of $H_n(F)$, denote respectively by 
\[
L_{\mathfrak{D}} :=(L(\Delta_i,\,\Delta_j))_{i,j=1}^\mu, \quad 
I_{\mathfrak{D}} := (\langle \Delta_i,\,\Delta_j\rangle)_{i,j=1}^\mu, 
\]
the matrices of the Seifert form and of the intersection form $\langle\cdot,\cdot \rangle$. 
Then, $L$ satisfies

\begin{prop}
Denote by $T_{\mathfrak{D}}$ the monodromy operator matrix. 
\begin{itemize}
\item[$\bullet$] $\det{L_{\mathfrak{D}}}=\pm1$, and $L_{\mathfrak{D}}$ is upper triangular. 
\item[$\bullet$] $I_{\mathfrak{D}}=-L_{\mathfrak{D}} -(-1)^n\cdot{}^tL_{\mathfrak{D}}$, and $T_{\mathfrak{D}}=(-1)^nL_{\mathfrak{D}}^{-1}{}^tL_{\mathfrak{D}}$. \QED
\end{itemize}
\end{prop}

\subsection{The Seifert form of a quasi-homogeneous IHS}

We summarize results in~\cite{S77-2} and~\cite{Sae00} for later use. 
Consider an IHS $(f=0)$. 

For a quasi-homogeneous polynomial $f$ with isolated singularity at $0$, we can define the {\it Jacobian ideal} of $f$ by 
\[
\mathrm{Jac}(f):= \left( \frac{\partial f}{\partial z_0},\cdots, \frac{\partial f}{\partial z_n}\right). 
\]
Denote by $J:=\mathbb{C}[z_0,\ldots,z_n]\slash\mathrm{Jac}(f)$ the {\it Jacobi algebra} graded by the weight system $\mbi{w}$, and $J_m$ the degree-$m$ part of $J$, for $m\geq0$. 

\begin{defn}[Poincar\'e series]
The Poincar\'e series of a graded algebra $A=\bigoplus_{n\in\mathbb{Q}} A_n$ is a power series 
\[
p_A(t) = \sum_{n\in\mathbb{Q}} (\dim{A_n})\,t^n. \quad \blacksquare
\]
\end{defn}

The Poincar\'e series of the Jacobi algebra for an IHS $X$ with weight system $\mbi{w}$ is given in ~\cite[Example (5.11)]{S77-2}: 
\begin{equation}
p_J(t) = \prod_{i=0}^n \frac{t^{w_i} - t}{1-t^{w_i}}. \label{PoincareSeries}
\end{equation}



Recall the formula for the $\mathbb{R}$-extension of the Seifert form $L_{\mathbb{R}}:=L\otimes_{\mathbb{Z}}\mathbb{R}$ given in~\cite{Sae00}. 


\begin{prop}[{\cite[Proposition 4.2]{Sae00}}]\label{Sae4.2}
The real Seifert form $L_{\mathbb{R}}$ is equivalent to the matrix 
\[
\bigoplus_{c_{\alpha}\not=0, \, \alpha\leq(n+1)/2} L^n_{\mathbb{R}} (c_{\alpha},\,\alpha), \quad {\it where}
\]
if $n$ is even, then, 
\[
L^n_{\mathbb{R}} (c_{\alpha},\,\alpha) := 
\begin{cases}
\bigoplus_{c_\alpha}
\begin{tiny}
\begin{pmatrix} 
\sin{\pi\alpha} & \cos{\pi\alpha} \\
-\cos{\pi\alpha} & \sin{\pi\alpha}
\end{pmatrix} \end{tiny}
& \textnormal{if } \alpha < \frac{n+1}{2}, \\
\bigoplus_{c_\alpha}(1) & \textnormal{if } \alpha = \frac{n+1}{2}, \textnormal{and } n\equiv0\mod4, \\
\bigoplus_{c_\alpha}(-1) & \textnormal{if } \alpha = \frac{n+1}{2}, \textnormal{and } n\equiv2\mod4, \\
\end{cases}
\]
if $n$ is odd, then, 
\[
L^n_{\mathbb{R}} (c_{\alpha},\,\alpha) := 
\begin{cases}
\bigoplus_{c_\alpha}
\begin{tiny}
\begin{pmatrix} 
\cos{\pi\alpha} & -\sin{\pi\alpha} \\
\sin{\pi\alpha} & \cos{\pi\alpha}
\end{pmatrix} \end{tiny} 
& \textnormal{if } \alpha < \frac{n+1}{2}, \\
\bigoplus_{c_\alpha}(1) & \textnormal{if } \alpha = \frac{n+1}{2}, \textnormal{and } n\equiv3\mod4, \\
\bigoplus_{c_\alpha}(-1) & \textnormal{if } \alpha = \frac{n+1}{2}, \textnormal{and } n\equiv1\mod4. \textnormal{\quad \QED}
\end{cases}
\]
\end{prop}


\begin{remark}
Let $L_\alpha$ be the eigenspace of eigenvalue $\alpha$ of the real Seifert form of a simple $K3$ singularity. 
It is easily seen by Propotision~\ref{Sae4.2} and basic linear algebra that 
\begin{itemize}
\item[$\bullet$] $c_2:=\dim{L_1}$ is none other than the coefficient of $t^2$ in $p_J(t)$. 
\item[$\bullet$] $\dim{L_{-1}}=2c_1$, where $c_1$ is the coefficient of $t$ in $p_J(t)$. 
In fact, we have $c_1=1$. 
\end{itemize}
\end{remark}

\begin{eg}\label{EgK3}
Consider a simple $K3$ singularity defined by $f=x^3 + y^4 + z^4 + w^6$ with weight system $\mbi{w}=(1/3,\, 1/4,\, 1/4,\, 1/6)$, and $Q=(4,3,3,2)$. 
The Poincar\'e series is 
\begin{tiny}
\begin{eqnarray*}
p_J(t) & = & \frac{t^{1/3}-t}{1-t^{1/3}}\cdot\frac{t^{1/4}-t}{1-t^{1/4}}\cdot\frac{t^{1/4}-t}{1-t^{1/4}}\cdot\frac{t^{1/6}-t}{1-t^{1/6}} \\
 & = & t + t^{7/6} + 2 t^{5/4} + 2 t^{4/3} + 2 t^{17/12} + 5 t^{3/2} + 4 t^{19/12} + 5 t^{5/3} \\
 & & + 6 t^{7/4} + 6 t^{23/12} + 8 t^2 + 6 t^{25/12} + 7 t^{13/6} + 6 t^{9/4} + 5 t^{7/3} \\
 & & + 4 t^{29/12} + 5 t^{5/2} + 2 t^{31/12} + 2 t^{8/3} + 2 t^{11/4} + t^{17/6} + t^3.  
\end{eqnarray*}
\end{tiny}
By Proposition \ref{Sae4.2}, the matrix $L_\mathbb{R}$ of the real Seifert form is equivalent to
\begin{tiny}
\begin{eqnarray*}
L_\mathbb{R} & \simeq & \bigoplus_{i=1}^8 (1)\oplus \begin{pmatrix} -1 & 0 \\ 0 & -1\end{pmatrix} \oplus 
\frac{1}{2}\begin{pmatrix} -\sqrt{2} & 1 \\ -1 & -\sqrt{2}\end{pmatrix} \oplus
\bigoplus_{i=1}^2 \frac{1}{\sqrt{2}} \begin{pmatrix} -1 & 1 \\ -1 & -1  \end{pmatrix} \oplus
\bigoplus_{i=1}^2 \frac{1}{\sqrt{2}} \begin{pmatrix} -1 & \sqrt{3} \\ -\sqrt{3} & -1\end{pmatrix} \oplus 
\\  & & 
\bigoplus_{i=1}^2 \frac{1}{2\sqrt{2}} \begin{pmatrix} 1-\sqrt{3} & 1+\sqrt{3} \\ -(1+\sqrt{3}) & 1-\sqrt{3}\end{pmatrix} \oplus 
\bigoplus_{i=1}^5 \begin{pmatrix}0 & 1 \\ -1 & 0 \end{pmatrix}  \oplus 
\bigoplus_{i=1}^4 \frac{1}{2\sqrt{2}}\begin{pmatrix} -1+\sqrt{3} & 1 + \sqrt{3} \\ -1-\sqrt{3} & -1+\sqrt{3} \end{pmatrix} \oplus 
\\ & & 
\bigoplus_{i=1}^5 \frac{1}{2}\begin{pmatrix} 1 & \sqrt{3} \\ -\sqrt{3} & 1 \end{pmatrix} \oplus 
\bigoplus_{i=1}^6 \frac{1}{\sqrt{2}}\begin{pmatrix} 1 & 1 \\ -1 & 1 \end{pmatrix} \oplus 
\bigoplus_{i=1}^7 \frac{1}{2}\begin{pmatrix} \sqrt{3} & 1 \\ -1 & \sqrt{3} \end{pmatrix} \oplus
\bigoplus_{i=1}^6 \frac{1}{2\sqrt{2}}\begin{pmatrix} 1+\sqrt{3} & -1+\sqrt{3} \\ 1-\sqrt{3} & 1 + \sqrt{3} \end{pmatrix}.
\end{eqnarray*}
\end{tiny}
The matrix $L_\mathbb{R}$ is of rank $90$ and determinant $1$, and has eigenvalues $-1$ and $1$, of which the eigenspaces are respectively of dimension 2, and 8. 
On the other hand, the transcendental lattice of the family has signature $(2,\, 8)$. 

\end{eg}

\subsection{Main Result}
Let $c_2$ be the dimension of the eigenspace of eigenvalue 1 of the real Seifert form of a simple $K3$ singularity, and $\rho$ the Picard number of the family of associated weighted $K3$ surfaces. 
Denote by $l(\Delta^{[1]})$ the total number of lattice points on all edges of the Newton polytope of the defining polynomial of the general section in the family. 
Generalizing Example \ref{EgK3}, we have our main result: 
\begin{prop}\label{MainThmSeifert}
The following formula holds for all simple $K3$ hypersurface singularities. 
\[
c_2 = l(\Delta^{[1]})-3 = 20 - \rho. 
\]
\end{prop}
\noindent
\proof
The value $c_2$ can be seen by finding the Poincar\'e series by formula~(\ref{PoincareSeries})\footnote{The full list of Poincar\'e series is given in Table \ref{PoiChar}. }. 
The Picard number $\rho$ is given in~\cite{Be02}. 
The value $l(\Delta^{[1]})$ is obtained by looking at the Newton polytope of the family~\cite{MaseDis}. 
Therefore, the statements follow by a direct computation. 
\QED

\setlength{\oddsidemargin}{-7mm}
\setlength{\textwidth}{17cm}
\setlength{\extrarowheight}{3.5pt}
\begin{landscape}
\tiny
\begingroup
\captionof{table}{{The Poincar\'e series and characteristic polynomials}}\label{PoiChar}
\endgroup
\begin{longtable}{llp{160mm}ccc}
\textnormal{No.} & $\mbi{w}$, $Q$ and $f$ & $p_J(t)$ and $\Div\Delta$ & $c_2$ & $\sign{(T)}$ & $l(\Delta^{[1]})$ \\ 
\hline
1. & $\f1{4}, \f1{4}, \f1{4}, \f1{4}$; $1,1,1,1$ & $10 t^{3/2}+10 t^{5/2}+4 t^{5/4}+16 t^{7/4}+16 t^{9/4}+4 t^{11/4}+t^3+19 t^2+t$ & $19$ & $(2, 19)$ & 22 \\ 
\cline{3-3} 
 & $w^4 + x^4 + y^4 + z^4$ & $20\Lambda_{4} + \Lambda_1$ \\ 
\hline

2. & $\f1{6}, \f1{4}, \f1{4}, \f1{3}$; $2,3,3,4$ & $5 t^{3/2}+5 t^{5/2}+2 t^{4/3}+5 t^{5/3}+5 t^{7/3}+2 t^{8/3}+2 t^{5/4}+6 t^{7/4}+6 t^{9/4}+2 t^{11/4}+t^{7/6}+7 t^{11/6}+7 t^{13/6}+t^{17/6}+2 t^{17/12}+4 t^{19/12}+6 t^{23/12}+6 t^{25/12}+4 t^{29/12}+2 t^{31/12}+t^3+8 t^2+t$ & $8$ & $(2, 8)$ & 11 \\ 
\cline{3-3} 
 & $w^6 + x^4 + y^4 + z^3$ & $6\Lambda_{12}+2\Lambda_6+2\Lambda_4-\Lambda_3+\Lambda_1$ \\ 
\hline

3. & $\f1{6}, \f1{6}, \f1{3}, \f1{3}$; $1,1,2,2$ & $8 t^{3/2}+8 t^{5/2}+5 t^{4/3}+12 t^{5/3}+12 t^{7/3}+5 t^{8/3}+2 t^{7/6}+14 t^{11/6}+14 t^{13/6}+2 t^{17/6}+t^3+16 t^2+t$ & $16$ & $(2, 16)$ & 19 \\
\cline{3-3} 
 & $w^6 + x^6 + y^3 + z^3$ & $16\Lambda_{6} +\Lambda_3 + \Lambda_1$ \\ 
\hline 

4. & $\f1{12}, \f1{4}, \f1{3}, \f1{3}$; $1,3,4,4$ & $5 t^{3/2}+5 t^{5/2}+4 t^{4/3}+8 t^{5/3}+8 t^{7/3}+4 t^{8/3}+2 t^{5/4}+8 t^{7/4}+8 t^{9/4}+2 t^{11/4}+t^{7/6}+10 t^{11/6}+10 t^{13/6}+t^{17/6}+t^{13/12}+4 t^{17/12}+7 t^{19/12}+10 t^{23/12}+10 t^{25/12}+7 t^{29/12}+4 t^{31/12}+t^{35/12}+t^3+10 t^2+t$ & $10$ & $(2, 10)$ & 13  \\
\hline
\end{longtable}
\end{landscape}

\setlength{\baselineskip}{16pt}

Makiko Mase, {\bf mtmase@arion.ocn.ne.jp}, Tokyo Metropolitan University. \\
1-1 Hachioji-shi Minami-osawa Tokyo, Japan 192-0397. 

\end{document}